
\input amstex 
\input amsppt.sty
\hsize 30pc
\vsize 47pc
\def\nmb#1#2{#2}         
\def\cit#1#2{\ifx#1!\cite{#2}\else#2\fi} 
\def\totoc{}             
\def\idx{}               
\def\ign#1{}             

\redefine\o{\circ}

\define\al{\alpha}

\define\ze{\zeta}

\define\io{\iota}

\define\la{\lambda}

\define\si{\sigma}

\define\ph{\varphi}

\define\ps{\psi}
\define\om{\omega}

\define\De{\Delta}

\define\Om{\Omega}
\predefine\ii\i
\redefine\i{^{-1}}
\define\row#1#2#3{#1_{#2},\ldots,#1_{#3}}
\define\x{\times}

\define\Der{{\operatorname{Der}}}
\define\Hom{\operatorname{Hom}}
\define\sign{\operatorname{sign}}

\define\End{\operatorname{End_{\Bbb K}}}
\define\ad{\operatorname{ad}}

\define\Lip{\operatorname{Lip}}
\define\Int{\operatorname{Int}}
\define\Out{{\operatorname{Out}}}
\define\Diag{{\operatorname{Diag}}}
\redefine\L{{\Cal L}}
\def\today{\ifcase\month\or
 January\or February\or March\or April\or May\or June\or
 July\or August\or September\or October\or November\or December\fi
 \space\number\day, \number\year}
\hyphenation{ho-mo-mor-phism boun-ded}
\topmatter
\title
The Fr\"olicher-Nijenhuis bracket for derivation based
non commutative differential forms \endtitle
\author  Michel Dubois-Violette\\
Peter W. Michor  \endauthor
\leftheadtext{\smc M\. Dubois-Violette, P\. Michor}
\rightheadtext{\smc Non-commutative Fr\"olicher-Nijenhuis bracket}
\affil
LPTHE Universit\'e Paris XI, B\^atiment 211, F-91405 Orsay Cedex,
France\\
Erwin Schr\"odinger International Institute of Mathematical Physics,
Wien, Austria
\endaffil
\address
M\. Dubois-Violette:
Laboratoire de Physique Th\'eorique et Hautes Energies,
Universit\'e Paris XI, B\^atiment 211, F-91405 Orsay Cedex,
France
\endaddress
\email flad\@qcd.th.u-psud.fr \endemail
\address
P\. Michor: Institut f\"ur Mathematik, Universit\"at Wien,
Strudlhofgasse 4, A-1090 Wien, Austria; and:
Erwin Schr\"odinger International Institute of Mathematical Physics,
Pasteurgasse 6/7, A-1090 Wien, Austria
\endaddress
\email Peter.Michor\@esi.ac.at \endemail
\date {\today} \enddate
\thanks
The work of both authors was supported by the program
`Geometry and Gravity', at the Isaac Newton Institute for
Mathematical Sciences, Cambridge, UK
\endthanks
\abstract
In commutative differential geometry the
Fr\"olicher-Nijenhuis bracket computes all kinds of curvatures and
obstructions to integrability. In \cit!{3} the Fr\"olicher-Nijenhuis
bracket was developped for universal differential forms of
non-commutative algebras, and several applications were given.
In this paper this bracket and the Fr\"olicher-Nijenhuis calculus will be
developped for several kinds of differential graded algebras based on
derivations, which were indroduced by \cit!{6}.
\endabstract

\endtopmatter

\document

\heading Table of contents \endheading
\noindent 1. Introduction \leaders \hbox to
     1em{\hss .\hss }\hfill {\eightrm 2}\par
\noindent 2. Convenient vector spaces \leaders \hbox to
     1em{\hss .\hss }\hfill {\eightrm 3}\par
\noindent 3. Preliminaries: graded differential algebras, \par
     derivations, and operations of Lie algebras \leaders \hbox to
     1em{\hss .\hss }\hfill {\eightrm 7}\par
\noindent 4. Derivations on universal differential forms
     \leaders \hbox to 1em{\hss .\hss }\hfill {\eightrm 9}\par
\noindent 5. The Fr{\accent "7F o}licher-Nijenhuis calculus on
     Chevalley type cochains \leaders \hbox to
     1em{\hss .\hss }\hfill {\eightrm 12}\par
\noindent 6. Description of all derivations \par
     in the Chevalley differential graded algebra \leaders \hbox to
     1em{\hss .\hss }\hfill {\eightrm 18}\par
\noindent 7. Diagonal bimodules \leaders \hbox to
     1em{\hss .\hss }\hfill {\eightrm 20}\par
\noindent 8. Derivations on the differential graded algebra
     $\Omega _{\Der}(A)$ \leaders \hbox to
     1em{\hss .\hss }\hfill {\eightrm 23}\par
\noindent 9. The differential graded algebra
     $\Omega _{\Out}(A)$ \leaders \hbox to
     1em{\hss .\hss }\hfill {\eightrm 27}\par

\head\totoc\nmb0{1}. Introduction \endhead

There are several generalizations of the differential calculus of
differential forms in the non-commutative setting \cit!{4},
\cit!{12}, \cit!{13}, \cit!{6}, \cit!{3}. We concentrate here on the
differential calculus based on derivations as generalizations of
vector fields, \cit!{6}.

Let us recall what are the the relevant notions of differential forms
in this context. Let $A$ be an associative algebra with a unit $1$.
As usual we think of $A$ as a generalization of an algebra of smooth
functions. Then it is natural to consider the Lie algebra $\Der(A)$
of all derivations of $A$ with values in $A$, i\.e\. of all
infinitesimal automorphisms of $A$, as a generalization of the Lie
algebra of vector fields. Alternatively, for an `invariant' theory
one may prefer to use the Lie algebra $\Out(A)$ of all derivations of
$A$ modulo the inner derivations. Then $\Out(A)$ is Morita invariant
and it also coincides with the Lie algebra of all vector fields in
case that $A$ is the algebra of smooth functions on a manifold.
As for the commutative case (see \cit!{14}) the notions of
differential forms can be extracted from the differential algebra
$C(\Der(A),A)$ of Chevalley-Eilenberg cochains of the Lie algebra
$\Der(A)$ with values in the $\Der(A)$-module $A$. There are then two
natural generalizations of the graded differenial algebra of
differential forms: A minimal one, $\Om_\Der(A)$, which is the
smallest differential subalgebra of $C(\Der(A),A)$ which contains
$A=C^0(\Der(A),A)$. And a maximal one,
$\underline\Om_\Der(A) =C_{Z(A)}(\Der(A),A)$,
which consists of all cochains in $C(\Der(A),A)$ which are module
homomorphisms for the module structure of $\Der(A)$ over the center
$Z(A)$ of $A$. In order to pass to the corresponding notions for
$\Out(A)$ we notice that there is a canonical operations, in the
sense of H\. Cartan \cit!{1}, \cit!{2}, $X\mapsto i_X$ for
$X\in \Der(A)$, of the Lie algebra $\Der(A)$ in the graded
differential algebra $C(\Der(A),A)$. Both $\Om_\Der(A)$ and
$\underline \Om_\Der(A)=C_{Z(A)}(\Der(A),A)$ are stable under this
operation and we define $\Om_\Out(A)$ and $\underline\Om_\Out(A)$ to
be the differential subalgebras of $\Om_\Der(A)$ and
$\underline\Om_\Der(A)$ consisting of all elements which are basic
with respect to the corresponding operation of the ideal $\Int(A)$ of
inner derivations of $\Der(A)$; one defines similarly
$C_{\Out(A)}(\Der(A),A)$. These graded differential algebras are also
obvious generalizations of differential forms. Notice, however, that
in contrast to $\Om_\Der(A)$ and $\underline\Om_\Der(A)$ there is in
general no differential calculus starting with $A$ in these algebras,
since they do not contain $A$ but merely its center $Z(A) =
\Om^0_\Out(A) = \underline\Om^0_\Out(A) = C^0_{\Out(A)}(\Der(A),A)$.
The differential algebra $C_{\Out(A)}(\Der(A),A)$ turns out out to be
the differentail algebra $C(\Out(A),Z(A))$ of all cochains of the Lie
algebra $\Out(A)$ with values in the center $Z(A)$, which is an
$\Out(A)$-module. Under this identification $\underline\Om_\Out(A)$
becomes the differential algebra $C_{Z(A)}(\Out(A),Z(A))$ of
$Z(A)$-multilinear cochains. This implies that
$\underline\Om_\Out(A)$ is a Morita invariant generalization of the
differential algebra of differential forms. $C(\Out(A),Z(A))$ is of
course also Morita invariant.

We shall develop the theory in several directions. Firstly, we shall
show that the derivation $d:A\to \Om^1_\Der(A)$ is universal for the
derivations of $A$ in a category of bimodules containing all
bimodules  which are isomorphic to sub bimodules of arbitrary
products of $A$  considered as a bimodule.: As suggested by
A\.~Connes, we call these last bimodules \idx{\it diagonal
bimodules}. This means that when one restricts attention to the above
catory of bimodules  containing the diagonal ones, the universal
property of the universal derivation $d:A\to \Om^1(A)$  factors
through the canonical surjective bimodule homomorphism
$\ze:\Om^1(A)\to \Om^1_\Der(A)$.

Secondly, we shall generalize for these differential forms the
Fr\"olicher-Nijenhuis calculus for vector valued differential forms.
As generalization of the space of vector valued differential forms we shall
consider $\Der(A,\Om_\Der(A))$ in the case of $\Om_\Der(A)$, and
$\Der(A,\underline\Om_\Der(A))$ in the case of
$\underline\Om_\Der(A)$, etc. For the universal differential
envelopping algebra $\Om(A)$ of $A$, the generalization of the
Fr\"olicher-Nijenhuis bracket has already been introduced on
$\Der(A,\Om_(A))$ in \cit!{3}, and we shall also define such a
generalization for $\Der(A,C(\Der(A),A))\cong C(\Der(A),\Der(A))$.
The generalizations proposed are natural, so that under the sequence
of homomorphisms and inclusions of graded differential algebras
$$
\Om(A)\to \Om_\Der(A) \subset \underline\Om_\Der(A) \subset C(\Der(A),A)
$$
the corresponding sequence
$$
\Der(A,\Om(A))\to \Der(A,\Om_\Der(A))
\subset \Der(A,\underline\Om_\Der(A)) \subset \Der(A,C(\Der(A),A))
$$
is a sequence of homomorphisms for the generalized
Fr\"olicher-Nijenhuis brackets.
Moreover we present here a novel approach to the
Fr\"olicher-Nijenhuis bracket, which uses the Chevalley coboundary
operator for the adjoint representation and which works also in other
situations, see the proofs of \nmb!{5.7} and \nmb!{5.8}.

Since it is useful to have a theory which is well suited to
topological algebras we develop from the beginning the whole theory
in the setting of convenient vector spaces as developped by
Fr\"olicher and Kriegl. The reasons for this are the following: If
the non-commutative theory should contain some version of
differential geometry, a manifold $M$ should be represented by the
algebra $C^\infty(M,\Bbb R)$ of smooth functions on it. The simplest
considerations of groups need products, and $C^\infty(M\x N,\Bbb R)$ is a
certain completion of the algebraic tensor product
$C^\infty(M,\Bbb R)\otimes C^\infty(N,\Bbb R)$. Now the setting of
convenient vector spaces offers in its multilinear version a
monoidally closed category, i\.e\. there is an appropriate tensor
product which has all the usual (algebraic) properties with respect
to bounded multilinear mappings. So multilinear algebra is carried
into this kind of functional analysis without loss. The theory of
convenient vector spaces is sketched in the first section \nmb!{2}.

We note that all results of this paper also hold in a purely
algebraic setting: Just equip each vector space with the finest
locally convex topology, then all linear mappings are bounded.

\heading\totoc\nmb0{2}. Convenient vector spaces \endheading

\subheading{\nmb.{2.1}} The traditional differential calculus works
well for Banach spaces. For more general locally convex spaces a
whole flock of different theories were developed, each of them rather
complicated and none really convincing. The main difficulty is that
the composition of linear mappings stops to be jointly continuous at
the level of Banach spaces, for any compatible topology. This was the
original motivation for the development of a whole new field within
general topology, convergence spaces.

Then in 1982, Alfred Fr\"olicher and Andreas Kriegl presented
independently the solution to the quest for the right differential
calculus in infinite dimensions. They joined forces in the further
development of the theory and the (up to now) final outcome is the
book \cit!{9}. See also the forthcoming book
\cit!{17}, or \cit!{15}, \cit!{16} for the material presented in this
section.

The appropriate spaces for this differential calculus are the
convenient vector spaces mentioned above. In addition to their
importance for differential calculus these spaces form a category
with very nice properties.

In this section we will sketch the basic definitions and the most
important results concerning convenient vector spaces and
Fr\"olicher-Kriegl calculus. All locally convex spaces will be
assumed to be Hausdorff.

\subheading{\nmb.{2.2}. The $c^\infty$-topology} Let $E$ be a
locally convex vector space. A curve $c:\Bbb R\to E$ is called
{\it smooth} or $C^\infty$ if all derivatives exist (and are
continuous) - this is a concept without problems. Let
$C^\infty(\Bbb R,E)$ be the space of smooth curves. It can be
shown that $C^\infty(\Bbb R,E)$ does not depend on the locally convex
topology of $E$, only on its associated bornology (system of bounded
sets).

The final topologies with respect to the following sets of mappings
into E coincide:
\roster
\item $C^\infty(\Bbb R,E)$.
\item Lipschitz curves (so that $\{\frac{c(t)-c(s)}{t-s}:t\neq s\}$
     is bounded in $E$).
\item $\{E_B\to E: B\text{ bounded absolutely convex in }E\}$, where
     $E_B$ is the linear span of $B$ equipped with the Minkowski
     functional $p_B(x):= \inf\{\la>0:x\in\la B\}$.
\item Mackey-convergent sequences $x_n\to x$ (there exists a sequence
     $0<\la_n\nearrow\infty$ with $\la_n(x_n-x)$ bounded).
\endroster
This topology is called the $c^\infty$-topology on $E$ and we write
$c^\infty E$ for the resulting topological space. In general (on the
space $\Cal D$ of test functions for example) it is finer than the
given locally convex topology; it is not a vector space topology,
since addition is no longer jointly continuous. The finest among all
locally convex topologies on $E$ which are coarser than the
$c^\infty$-topology is the bornologification of the given locally
convex topology. If $E$ is a Fr\'echet space, then $c^\infty E = E$.

\subheading{\nmb.{2.3}. Convenient vector spaces} Let $E$ be a
locally convex vector space. $E$ is said to be a {\it convenient
vector space} if one of the following equivalent
conditions is satisfied (called $c^\infty$-completeness):
\roster
\item Any Mackey-Cauchy-sequence (so that $(x_n-x_m)$ is Mackey
     convergent to 0) converges.
\item If $B$ is bounded closed absolutely convex, then $E_B$ is a
     Banach space.
\item Any Lipschitz curve in $E$ is locally Riemann integrable.
\item For any $c_1\in C^\infty(\Bbb R,E)$ there is
     $c_2\in C^\infty(\Bbb R,E)$ with $c_1=c_2'$ (existence of
     antiderivative).
\endroster

Obviously $c^\infty$-completeness is weaker than
sequential completeness so any  sequentially complete locally convex
vector space is convenient.
 From \nmb!{2.2}.4 one easily sees that $c^\infty$-closed linear
subspaces of convenient vector spaces are again convenient. We
always assume that a convenient vector space is equipped with its
bornological topology.

\proclaim{\nmb.{2.4}. Lemma} Let $E$ be a locally convex space.
Then the following properties are equivalent:
\roster
\item $E$ is $c^\infty$-complete.
\item If $f:\Bbb R\to E$ is scalarwise $\Lip^k$, then $f$ is
     $\Lip^k$, for $k>1$.
\item If $f:\Bbb R\to E$ is scalarwise $C^\infty$ then $f$ is
     differentiable at 0.
\item If $f:\Bbb R\to E$ is scalarwise $C^\infty$ then $f$ is
     $C^\infty$.
\endroster
\endproclaim
Here a mapping $f:\Bbb R\to E$ is called $\Lip^k$ if all partial
derivatives up to order $k$ exist and are Lipschitz, locally on
$\Bbb R$. $f$ scalarwise $C^\infty$ means that $\la\o f$ is $C^\infty$
for all continuous linear functionals on $E$.

This lemma says that on a convenient vector space one can recognize
smooth curves by investigating compositions with continuous linear
functionals.

\subheading{\nmb.{2.5}. Smooth mappings} Let $E$ and $F$ be locally
convex vector spaces. A mapping $f:E\to F$ is called {\it smooth} or
$C^\infty$, if $f\o c\in C^\infty(\Bbb R,F)$ for all
$c\in C^\infty(\Bbb R,E)$; so
$f_*: C^\infty(\Bbb R,E)\to C^\infty(\Bbb R,F)$ makes sense.
Let $C^\infty(E,F)$ denote the space of all smooth mappings from $E$
to $F$.

For $E$ and $F$ finite dimensional this gives the usual notion of
smooth mappings.
Constant mappings are smooth. Multilinear mappings are smooth if and
only if they are bounded. Therefore we denote by $L(E,F)$ the space
of all bounded linear mappings from $E$ to $F$.

\proclaim{\nmb.{2.6}. Lemma }
For any locally convex space $E$ there is a convenient vector space
$\tilde E$ called the completion of $E$ and a bornological embedding
$i:E\to \tilde E$, which is characterized by the
property that any bounded linear map from $E$ into an arbitrary
convenient vector space extends to $\tilde E$.
\endproclaim

\subheading{\nmb.{2.7} }
As we will need it later on we describe the completion in a special
situation:
Let $E$ be a locally convex space with completion $i:E\to \tilde E$,
$f:E\to E$ a bounded projection and $\tilde f:\tilde E\to \tilde E$ the
prolongation of $i\o f$. Then $\tilde f$ is also a projection and
$\tilde f(\tilde E)=\ker (Id-\tilde f)$ is a $c^\infty$-closed and
thus convenient linear subspace of $\tilde E$. Using that $f(E)$ is a
direct summand in $E$ one easily shows that $\tilde f(\tilde E)$ is the
completion of $f(E)$. This argument applied to $Id-f$ shows that
$\ker (\tilde f)$ is the completion of $\ker (f)$.

\subheading{\nmb.{2.8}. Structure on $C^\infty(E,F)$} We equip the
space $C^\infty(\Bbb R,E)$ with the bornologification of the topology
of uniform convergence on compact sets, in all derivatives
separately. Then we equip the space $C^\infty(E,F)$ with the
bornologification of the initial topology with respect to all
mappings $c^*:C^\infty(E,F)\to C^\infty(\Bbb R,F)$, $c^*(f):=f\o c$,
for all $c\in C^\infty(\Bbb R,E)$.

\proclaim{\nmb.{2.9}. Lemma } For locally convex spaces $E$ and $F$
we have:
\roster
\item If $F$ is convenient, then also $C^\infty(E,F)$ is convenient,
     for any $E$. The space $L(E,F)$ is a closed linear subspace of
     $C^\infty(E,F)$, so it is convenient also.
\item If $E$ is convenient, then a curve $c:\Bbb R\to L(E,F)$ is
     smooth if and only if $t\mapsto c(t)(x)$ is a smooth curve in $F$
     for all $x\in E$.
\endroster
\endproclaim

\proclaim{\nmb.{2.10}. Theorem} The category of convenient vector
spaces and smooth mappings is cartesian closed. So we have a natural
bijection
$$C^\infty(E\x F,G)\cong C^\infty(E,C^\infty(F,G)),$$
which is even a diffeomorphism.
\endproclaim

Of course this statement is also true for $c^\infty$-open subsets of
convenient vector spaces.

\proclaim{\nmb.{2.11}. Corollary } Let all spaces be convenient vector
spaces. Then the following canonical mappings are smooth.
$$\align
&\operatorname{ev}: C^\infty(E,F)\x E\to F,\quad
     \operatorname{ev}(f,x) = f(x).\\
&\operatorname{ins}: E\to C^\infty(F,E\x F),\quad
     \operatorname{ins}(x)(y) = (x,y).\\
&(\quad)^\wedge :C^\infty(E,C^\infty(F,G))\to C^\infty(E\x F,G),
	\quad \hat f(x,y)=f(x)(y).\\
&(\quad)\spcheck :C^\infty(E\x F,G)\to C^\infty(E,C^\infty(F,G)),
	\quad \check g(x)(y)=g(x,y).\\
&\operatorname{comp}:C^\infty(F,G)\x C^\infty(E,F)\to C^\infty(E,G)\\
&C^\infty(\quad,\quad):C^\infty(F,F')\x C^\infty(E',E)\to
     C^\infty(C^\infty(E,F),C^\infty(E',F'))\\
&\qquad (f,g)\mapsto(h\mapsto f\o h\o g)\\
&\prod:\prod C^\infty(E_i,F_i)\to C^\infty(\prod E_i,\prod F_i)
\endalign$$
\endproclaim

\proclaim{\nmb.{2.12}. Theorem} Let $E$ and $F$ be convenient vector
spaces. Then the differential operator
$$\gather d: C^\infty(E,F)\to C^\infty(E,L(E,F)), \\
df(x)v:=\lim_{t\to0}\frac{f(x+tv)-f(x)}t,
\endgather$$
exists and is linear and bounded (smooth). Also the chain rule holds:
$$d(f\o g)(x)v = df(g(x))dg(x)v.$$
\endproclaim

\subheading{\nmb.{2.13} }
The category of convenient vector spaces and bounded linear maps
is complete and cocomplete, so all categorical limits and colimits
can be formed. In particular we can form products and direct sums of
convenient vector spaces.

For convenient vector spaces $E_1$,\dots ,$E_n$ and $F$ we can now
consider the space of all bounded $n$-linear maps, $L(E_1,\dots ,E_n;F)$,
which is a closed linear subspace of $C^\infty(\prod _{i=1}^nE_i,F)$
and thus again convenient. It can be shown that multilinear maps are
bounded if and only if they are partially bounded, i\.e\. bounded in
each coordinate and that there is a natural isomorphism (of
convenient vector spaces) $L(E_1,\dots ,E_n;F)\cong
L(E_1,\dots ,E_k;L(E_{k+1},\dots ,E_n;F))$

\proclaim{\nmb.{2.14}. Theorem }
On the category of convenient vector spaces there is a unique tensor product
$\tilde \otimes$ which makes the category symmetric monoidally
closed, i\.e\. there are natural isomorphisms of convenient vector
spaces $L(E_1;L(E_2,E_3))\cong
L(E_1\tilde \otimes E_2,E_3)$,
$E_1\tilde \otimes E_2\cong E_2\tilde \otimes E_1$,
$E_1\tilde \otimes (E_2\tilde \otimes E_3)\cong (E_1\tilde \otimes
E_2)\tilde \otimes E_3$ and $E\tilde \otimes \Bbb R\cong E$.
\endproclaim
The tensor product can be constructed as follows: On the algebraic
tensor product put the finest locally convex topology such that the
canonical bilinear map from the product into the tensor product is
bounded and then take the completion of this space.

\subheading{\nmb.{2.15}. Remarks } Note that the conclusion of
theorem \nmb!{2.10} is the starting point of the classical calculus of
variations, where a smooth curve in a space of functions was assumed
to be just a smooth function in one variable more.

If one wants theorem \nmb!{2.10} to be true and assumes some other obvious
properties, then the calculus of smooth functions is already uniquely
determined.

There are, however, smooth mappings which are not continuous. This is
unavoidable and not so horrible as it might appear at first sight.
For example the evaluation $E\x E'\to\Bbb R$ is jointly continuous if
and only if $E$ is normable, but it is always smooth. Clearly smooth
mappings are continuous for the $c^\infty$-topology.
For Fr\'echet spaces smoothness in the sense described here coincides
with the usual notion of $C^\infty$.

\proclaim{\nmb.{2.16}. Lemma } {\rm \cit!{3}, 2.7}.
Let $A$ be a convenient algebra, $M$ a convenient right $A$-module
and $N$ a convenient left $A$-module. This means that all structure
mappings are bounded bilinear.
\roster
\item There is a convenient vector space $M\tilde \otimes _AN$ and a
     bounded bilinear map $b:M\times N\to M\tilde \otimes _AN$,
     $(m,n)\mapsto m\otimes _An$ such that
     $b(ma,n)=b(m,an)$ for all $a\in A$, $m\in M$ and $n\in N$ which has
     the following universal property: If $E$ is a convenient vector space
     and $f:M\times N\to E$ is a bounded bilinear map such that
     $f(ma,n)=f(m,an)$ then there is a unique bounded linear map
     $\tilde f:M\tilde \otimes _AN\to E$ with $\tilde f\o b=f$.
\item Let $L^A(M,N;E)$ denote the space of all bilinear bounded maps
     $f:M\times N\to E$ having the above property, which is a closed
     linear subspace of $L(M,N;E)$. Then we have an isomorphism of
     convenient vector spaces $L^A(M,N;E)\cong L(M\tilde \otimes _AN,E)$.
\item If $B$ is another convenient algebra such that $N$ is a
     convenient right $B$-module and such that the actions of $A$
     and $B$ on $N$ commute, then $M\tilde \otimes _AN$ is in a canonical
     way a convenient right $B$-module.
\item If in addition $P$ is a convenient left $B$-module then there
     is a natural isomorphism of convenient vector spaces
$$M\tilde \otimes _A(N\tilde \otimes _BP)\cong (M\tilde \otimes
     _AN)\tilde \otimes _BP$$
\endroster
\endproclaim

\subhead\nmb.{2.17}. Remark \endsubhead
In the following all spaces will be convenient spaces, and all
multilinear mappings will be bounded, even if not stated explicitly.
This includes the purely algebraic theory, where one just equips each
vector space with its finest locally convex topology, because then
each multilinear mapping is bounded automatically.

\head\totoc\nmb0{3}. Preliminaries: graded differential algebras, 
derivations, and operations of Lie algebras \endhead

\subhead\nmb.{3.1}. Graded derivations \endsubhead
Let $\frak A=\bigoplus_{k\in \Bbb Z}\frak A^k$ be a graded
associative algebra with unit, so that
$\frak A^k.\frak A^l\subset \frak A^{k+l}$.
We
denote by $\operatorname{Der}_k\frak A$ the space of
all \idx{\it (graded) derivations} of degree $k$, i\.e\. all linear
mappings $D:\frak A \to \frak A$ with $D(\frak A^l) \subset
\frak A^{k+l}$ and $D(\ph\wedge \ps) =
D(\ph)\wedge \ps +(-1)^{kl}\ph \wedge
D(\ps)$ for $\ph \in \frak A^l$.
Then the space
$\Der\frak A = \bigoplus_k\Der_k\frak A$ is a
graded Lie algebra with the graded commutator
$[D_1,D_2] := D_1\o D_2 - (-1)^{k_1k_2}D_2\o D_1$ as bracket.
This means that the bracket is graded anticommutative,
$[D_1,D_2] = -(-1)^{k_1k_2}[D_2,D_1]$, and satisfies the graded
Jacobi identity
$$[D_1,[D_2,D_3]] = [[D_1,D_2],D_3] + (-1)^{k_1k_2}[D_2,[D_1,D_3]]$$
(so that $ad(D_1) = [D_1,\quad]$ is itself a derivation of degree
$k_1$).

Let
$Z(\frak A)^q=\{a\in\frak A^q:
     [a,b]=ab-(-1)^{ql}ba=0\text{ for all }b\in \frak A^l\text{ and
     all }p\in \Bbb Z\}$
and consider the \idx{\it graded center}
$Z(\frak A)=\bigoplus_{q\in \Bbb Z}Z(\frak A)^q$ of $\frak A$, a
graded commutative algebra with unit, which is stable under
$\Der(\frak A)$.
Then $\Der(\frak A)$ is a (left) graded $Z(\frak A)$-module and we have
$[a.D_1,D_2] = a.[D_1,D_2] -(-1)^{(q+k_1)k_2}D_2(a).D_1$.

\subhead\nmb.{3.2}. Graded differential algebras \endsubhead
In this paper, a graded differential algebra is $\Bbb Z$-graded
associative algebra
$\frak A=\bigoplus_{k\in \Bbb Z} \frak A^k$
equipped with a graded derivation $d$ of degree $1$ with $d^2=0$. It
is called the \idx{\it differential} of $\frak A$.

\subhead\nmb.{3.3}. The graded differential algebra of
universal differential forms \endsubhead
Let $A$ be a convenient associative algebra with unit, and
let $(\Om^*(A),d)$ be the convenient graded differential algebra of
(universal) (K\"ahler) differential forms, see for example \cit!{3}.
Let us repeat quickly its construction: $\Om^1(A)$ is the kernel of
the multiplication $\mu:A\tilde\otimes A\to A$, a convenient
$A$-bimodule. The bounded linear mapping $d:A\to \Om^1(A)$, given by
$d(a)= 1\otimes a -a \otimes 1$, has the following universal
property, and the pair $(\Om^1(A),d)$ is uniquely determined by it:
\roster
\item For any bounded derivation $D:A\to N$ into a convenient $A$-bimodule
     $N$ there is a unique bounded $A$-bimodule homomorphism
     $j_D:\Om^1(A)\to N$ such that $D=j_D\o d$.
\endroster
We put $\Om^0(A):=A$, and for $k\in\Bbb Z$ we define
$\Om^k(A):=\Om^1(A)\tilde \otimes_A\dots\tilde \otimes_A\Om^1(A)$
($k$ factors).
There is a canonical extension of $d:A\to \Om^1(A)$ to a bounded
derivation of differential of graded algebra, i\.e\. a bounded
derivation degree 1 satisfying $d^2=0$.
$$
A @>d>> \Om^1(A) @>d>> \Om^2(A) @>d>>
\Om^3(A) @>d>> \dotso
$$
and the resulting convenient graded differential algebra has the
following universal property and is uniquely determined by it:
\roster
\item [2] For any bounded homomorphism $\ph:A\to B$ of convenient
     algebras and for any convenient graded differential algebra
     $(\Cal B = \bigoplus_{k=0}^{\infty}\Cal B^k,d^{\Cal B})$
     with  $\Cal B_0= B$ there exists a unique extension of $\ph$ to
     homomorphism $\Om(A)\to \Cal B$ of graded differential algebras.
\endroster

\subhead\nmb.{3.4}  \endsubhead
Let $A$ be a convenient associative algebra with unit, and denote its
center by $Z(A)$.
Let $\Der(A)=\{X\in L(A,A):X(ab)=X(a)b+aX(b)\text{ for all }a,b\in A\}$
be the convenient Lie algebra of all derivations of $A$.
Note that $\Der(A)$ is a bounded module over the center $Z(A)$ of $A$
only, with the usual formulas $[X,aY]=X(a).Y+ a[X,Y]$, etc.
Let $\frak g$ denote a convenient Lie algebra, which acts on $A$ by
derivations, i\.e\. there is a bounded Lie algebra homomorphism
$\frak g\to \Der(A)$.

We consider now the (Chevalley) graded differential algebra of
$\frak g$ with coefficients in $A$, which is given as follows.
Let $C^k(\frak g,A) := L^k_{\text{skew}}(\frak g,A)$
denote the convenient vector space of all bounded $k$-linear
(over $\Bbb R$) skew symmetric mappings $\frak g^k\to A$.
Let $C^0(\frak g,A):=A$. We put:
$$\align
C(\frak g,A) :&= \bigoplus_{k\ge0}C^k(\frak g,A)
     = \bigoplus_{k\ge0}L^k_{\text{skew}}(\frak g,A),\\
d\ph(X_0,\dots,X_k)
     :&= \sum_{i=0}^k(-1)^iX_i(\ph(X_0,\dots,\widehat{X_i},\dots,X_k))\\
&\quad + \sum_{i<j}(-1)^{i+j}\ph([X_i,X_j],X_0,\dots,\widehat{X_i},\dots,
        \widehat{X_j},\ldots,X_k),\\
(\ph.\ps)(X_1,\dots,X_{k+l})
     &= \tfrac 1{k!l !}\sum_{\si\in\Cal S_{k+l}}\sign\si.
\ph(X_{\si1},\dots,X_{\si k}).\ps(X_{\si(k+1)},\dotsc,X_{\si(k+l)}).
\endalign$$
Then $(C(\frak g,A),d)$ is a graded differential algebra.
It is a graded commutative differential algebra if and only if $A$ is
commutative.

\subhead\nmb.{3.5} \endsubhead
If $\frak g$ is a convenient Lie algebra and if $(\frak A,d)$ is a
graded differential algebra, an operation in the sense of Cartan of
$\frak g$ in $\frak A$ is a linear mapping
$\frak g\to \Der_{-1}(\frak A)$, written $X\mapsto i_X$, such that
for $\L_X:= i_X\,d + d\,i_X = [i_X,d]\in\Der_0(\frak A)$ we have for
all $X$, $Y\in \frak g$:
\roster
\item "" $i_X\,i_Y+i_Y\,i_X = [i_X,i_Y] = 0$
\item "" $\L_X,i_Y - i_Y\,\L_X = [\L_X,i_Y] = i_{[X,Y]}$.
\item "" $\L_X,\L_Y - \L_Y\,\L_X = [\L_X,\L_Y] = \L_{[X,Y]}$.
\endroster
An element $\al\in \frak A$ is called \idx{\it horizontal} with
respect to $\frak g$ whenever $i_X\al=0$ for all $X\in\frak g$; it is
called \idx{\it invariant} with respect to $\frak g$ if $\L_X\al=0$
for all $X\in\frak g$; and $\al$ is called \idx{\it basic} if it is
both horizontal and invariant.

The convenient space $\frak A_H$ of horizontal elements of $\frak A$ is a
graded
subalgebra of $\frak A$ which is stable under $\L_X$ for all
$X\in \frak g$. The convenient space $\frak A_I$ of all invariant elements is a
graded differential subalgebra of $(\frak A,d)$, and the convenient
space $\frak A_B$ of all basic elements is a graded
differential subalgebra of $\frak A_I$, thus also of $\frak A$.

\head\totoc\nmb0{4}. Derivations on universal differential forms \endhead

\subhead\nmb.{4.1} \endsubhead
Let $A$ be a convenient associative algebra with unit, and
let $(\Om^*(A),d)$ be the convenient graded differential algebra of
universal differential forms as described in \nmb!{3.3}.
In this section we review from \cit!{3} the description of all graded
derivations on the graded differential algebra $\Om^*(A)$.
This lead directly to what we like to call the `calculus of
Fr\"olicher-Nijenhuis'.

\subhead\nmb.{4.2} Derivations vanishing on $A$ \endsubhead
For every derivation $X\in \Der(A)$ there exists a bounded $A$-bimodule
homomorphism $j_X:\Om^1(A)\to A$, by the universal property
\nmb!{3.3}.(1). It prolongs uniquely to a graded derivation
$j(X)=j_X:\Om(A)\to\Om(A)$ of degree $-1$
which is called the
\idx{\it contraction operator} of the derivation $X$.

More generally, let us consider a graded derivation
$D\in\Der_k(\Om^*(A))$ of degree $k$ which vanishes on $A$,
$D\mid \Om^0(A) = 0$.
Such derivations are called \idx{\it algebraic derivations}; they
form a Lie subalgebra of $\Der(\Om^*(A))$.
Then $D(a\om) =
aD(\om)$ and $D(\om a)=D(\om)a$ for $a\in A$, so $D$ restricts to a
bounded bimodule homomorphism, an element of
$\Hom^A_A(\Om^l(A),\Om^{l+k}(A))$.
Since $\Om^l(A)$ for $l\ge 1$ is generated by $\Om^1(A)$,
the derivation $D$ is uniquely determined by its restriction
$D|\Om^1(A)\in\Hom^A_A(\Om^1(A),\Om^{k+1}(A))$.

Let us denote by
$$K:= (D|\Om^1(A))\o d \in \Der(A,\Om^{k+1}(A))$$
the corresponding derivation on $A$.
We write $D=j(K)=j_K$ to express the unique dependence of $D$ on $K$.
Note the defining equation $j_K(a_1.d(a_2).a_3)= a_1.K(a_2).a_3$ for
$a_i\in A$.

Conversely for $K \in \Der(A,\Om^{k+1}(A))$ with corresponding
homomorphism
$j_K\in \Hom^A_A(\Om^1(A),\Om^{k+1}(A))$ and $\om_i\in\Om^1(A)$
the formula
$$j_K(\om_0\otimes_A  \dots\otimes_A  \om_\ell)
     =\sum_{i=0}^\ell(-1)^{ik}\om_0\otimes_A  \dots
     \otimes_A  j_K(\om_i)\otimes_A  \dots \otimes_A \om_k$$
defines an algebraic graded derivation $j_K \in \Der_k\Om(A)$ and any
algebraic derivation is of this form.
The mapping
$$j:\Der(A,\Om^{k+1}(A))\to \Der_k\Om(A)$$
induces an isomorphism of convenient vector spaces onto the closed
linear subspace of $\Der_k\Om (A)$ consisting of all graded bounded
derivations which vanish on $A$.

By stipulating $j([K,L]^{\De}):= [j_K,j_L]$ we get a bracket
$[\quad,\quad]^{\De}$ on the space
$\Der(A,\Om^{*+1}(A))$ which defines a convenient
graded Lie algebra structure with the grading as indicated, and
for $K\in \Der(A,\Om^{k+1}(A))$, and
$L\in \Der(A,\Om^{k+1}(A))$ we have
$$[K,L]^{\De} = j_K\o L - (-1)^{k\ell}j_L\o K.$$
$[\quad,\quad]^{\De}$ is a version of the bracket of Gerstenhaber, De
Wilde - Lecomte, see \cit!{10}, \cit!{11}, and \cit!{5}.

\subhead\nmb.{4.3}. Lie derivations and the Fr\"olicher-Nijenhuis
bracket \endsubhead
The exterior derivative $d$ is an element of
$\Der_1\Om(A)$. We define for $K \in \Der(A,\Om^{k}(A))$
the \idx{\it Lie derivation}
$\L_K = \L(K) \in \Der_k\Om(A)$ by
$$\L_K := [j_K,d] = j_K\,d-(-1)^{k-1}d\,j_K.$$
Then the mapping $\L:\Der(A,\Om^{*}(A)) \to  \Der_*\Om(A)$
is obviously bounded
and it is injective by the universal property of $\Om^1(A)$, since
$\L_Ka=j_Kda = K(a)$ for $a\in A$. Note that $\L_K$ is an extension
of the derivation $K:A\to \Om^k(A)$ to a graded derivation of the
graded algebra $\Om^*(A)$ of degree $k$.

\proclaim{Lemma} {\rm \cit!{3}, 4.7.}
For any graded derivation $D \in \Der_k\Om(A)$
there are unique homomorphisms
$K \in \Der(A,\Om^{k}(A))$ and
$L \in \Der(A,\Om^{k+1}(A))$
such that $$D = \L_K + j_L.$$
We have $L=0$ if and only if $[D,d]=0$. $D$ is algebraic if and
only if $K=0$.\qed
\endproclaim

Note that $j_d\om = k\om$ for $\om \in \Om^k(A)$.
Therefore we have
$\L_d\om = j_d\,d\om - d\,j_d\om =
(k+1)d\om - kd\om = d\om$. Thus $\L_d = d$.

Let $K \in \Der(A,\Om^{k}(A))$ and
$L \in \Der(A,\Om^{l}(A))$.
Then obviously $[[\L_K,\L_L],d] =0$, so we have
$$[\L(K),\L(L)] = \L([K,L])$$
for a uniquely defined
$[K,L] \in \Der(A,\Om^{k+l}(A))$. This
vector valued form $[K,L]$ is called the \idx{\it (abstract)
Fr\"olicher-Nijenhuis bracket} of $K$ and $L$.

The space
$\Der(A,\Om^{*}(A)) = \bigoplus_k\Der(A,\Om^{k}(A))$
with its usual grading and the Fr\"olicher-Nijen\-huis
bracket is a convenient graded Lie algebra.
$d \in \Der(A,\Om^1(A))$ is in the center,
i.e. $[K,d] = 0$ for all $K$, see \cit!{3}, 4.9.

$\L:(\Der(A,\Om^*(A)), [\quad,\quad]) \to  \Der\Om(A)$ is
a bounded injective homomorphism of graded Lie algebras.
For $K \in \Der(A,\Om^k(A))$ and
$L \in \Der(A,\Om^{l+1}(A))$ we have
$$
[\L_K,j_L] = j([K,L]) - (-1)^{kl}\L(j_L\o K).\tag1
$$
For
$K_i \in \Der(A,\Om^{k_i}(A))$ and
$L_i \in \Der(A,\Om^{k_i+1}(A))$ we have
$$\align
&\left[\L_{K_1}+j_{L_1},\L_{K_2}+j_{L_2}\right] = \tag2 \\
&\qquad = \L\left([K_1,K_2]
    + j_{L_1}\o K_2 - (-1)^{k_1k_2}j_{L_2}\o K_1\right)\\
&\qquad\quad + i\left([L_1,L_2]^{\De}
    + [K_1,L_2] -(-1)^{k_1k_2}[K_2,L_1]\right).
\endalign$$
Each summand of this formula looks like a semidirect product of
graded Lie algebras, but the mappings
$$\align
j_*:\Der(A,\Om^{*-1}(A)) &\to
     \End(\Der(A,\Om^*(A)),[\quad,\quad]) \\
\ad:\Der(A,\Om^*(A)) &\to
     \End(\Der(A,\Om^{*-1}(A)),[\quad,\quad]^{\De}),\\
\ad_KL&=[K,L],
\endalign$$
do not take values in the subspaces of graded derivations. We
have instead for
$K \in \Der(A,\Om^k(A))$ and
$L \in \Der(A,\Om^{l+1}(A))$
the following relations:
$$\align &j_L\o [K_1,K_2] = [j_L\o K_1,K_2]
     + (-1)^{k_1l}[K_1,j_L\o K_2] \tag3 \\
&\qquad -\left((-1)^{k_1l}j(\ad_{K_1}L)\o K_2
      - (-1)^{(k_1+l)k_2}j(\ad_{K_2}L)\o K_1\right)\\ \allowdisplaybreak
&\ad_K[L_1,L_2]^{\De} = [\ad_KL_1,L_2]^{\De}
      + (-1)^{kk_1}[L_1,\ad_KL_2]^{\De} -\tag4 \\
&\qquad -\left((-1)^{kk_1}\ad(j(L_1)\o K)L_2
      - (-1)^{(k+k_1)k_2}\ad(j(L_2)\o K)L_1\right)
\endalign$$

\subheading{\nmb.{4.4}. Naturality of the Fr\"olicher-Nijenhuis
bracket} Let $f:A \to  B$ be a bounded algebra homomorphism.
Two elements
$K \in \Der(A,\Om^k(A))$ and
$K' \in \Der(B,\Om^k(B))$
are called \idx{\it $f$-related} or \idx{\it
$f$-dependent}, if we have
$$K'\o f = \Om^k(f)\o K: A\to\Om^k(B) ,\tag{1}$$
where $\Om^*(f):\Om^*(A)\to \Om^*(B)$ is given by the universal
property \nmb!{3.3}.\therosteritem2.
 From \cit!{3}, 4.12. we have the following results:
\roster
\item[2] If $K$ and $K'$ as above are $f$-related then
  $j_{K'}\o \Om(f) = \Om(f) \o j_{K}:\Om(A) \to  \Om(B)$.
\item If $j_{K'}\o \Om(f)| d(A) = \Om(f)\o j_{K}| d(A)$, then
  $K$ and $K'$ are $f$-related, where $d(A)\subset \Om^1(A)$
  denotes the space of exact $1$-forms.
\item If $K_j$ and $K_j'$ are $f$-related for $j=1,2$, then
  $j_{K_1}\o K_2$ and $j_{K'_1}\o K'_2$ are $f$-related, and also
  $[K_1,K_2]^{\De}$ and $[K'_1,K'_2]^{\De}$ are $f$-related.
\item If $K$ and $K'$ are $f$-related then
  $\L_{K'}\o \Om(f) = \Om(f)\o \L_{K}:\Om(A) \to  \Om(B)$.
\item If $\L_{K'}\o \Om(f)\mid \Om^0(A) = \Om(f)\o\L_{K}\mid \Om^0(A)$,
  then $K$ and $K'$ are $f$-related.
\item If $K_j$ and $K_j'$ are $f$-related for $j=1,2$, then
  their Fr\"olicher-Nijenhuis brackets $[K_1,K_2]$ and
  $[K'_1,K'_2]$ are also $f$-related.
\endroster

\head\totoc\nmb0{5}. The Fr\"olicher-Nijenhuis calculus
on Chevalley type cochains \endhead

\subhead\nmb.{5.1}  \endsubhead
Let $A$ be a convenient associative algebra with unit, with
center $Z(A)$, and with convenient Lie algebra of derivations
$\Der(A)=\{X\in L(A,A):X(ab)=X(a)b+aX(b)\text{ for all }a,b\in A\}$,
a bounded module over the center $Z(A)$ of $A$.
We consider now the (Chevalley) graded differential algebra of
$\Der(A)$ (see \nmb!{3.4}) which is given as follows.
Let $C^k(\Der(A),A) := L^k_{\text{skew}}(\Der(A),A)$
denote all bounded $k$-linear
(over $\Bbb R$) skew symmetric mappings $\Der(A)^k\to A$, and let
denote $C^k_{Z(A)}(\Der(A),A)$ denote the subspace of those mappings
which are $k$-linear over $Z(A)$. Let $C^0(\Der(A),A):=A$. We
put:
$$\align
C(\Der(A),A) :&= \bigoplus_{k\ge0}C^k(\Der(A),A)
     = \bigoplus_{k\ge0}L^k_{\text{skew}}(\Der(A),A),\\
d\ph(X_0,\dots,X_k)
     :&= \sum_{i=0}^k(-1)^iX_i(\ph(X_0,\dots,\widehat{X_i},\dots,X_k))\\
&\quad + \sum_{i<j}(-1)^{i+j}\ph([X_i,X_j],X_0,\dots,\widehat{X_i},\dots,
        \widehat{X_j},\ldots,X_k),\\
(\ph.\ps)(X_1,\dots,X_{k+l})
     &= \tfrac 1{k!l !}\sum_{\si\in\S_{k+l}}\sign\si.
\ph(X_{\si1},\dots,X_{\si k}).\ps(X_{\si(k+1)},\dotsc,X_{\si(k+l)}).
\endalign$$
Then $(C(\Der(A),A),d)$ is a graded differential algebra.
It is a graded commutative differential algebra if and only if $A$ is
commutative. If $A=C^\infty(M)$ for a smooth manifold $M$ then
$(C(\Der(A),A),d)$ is larger than the differential graded
algebra used for the Gel'fand-Fuks cohomology of the manifold $M$.

The subspace
$$
C_{Z(A)}(\Der(A),A) := \bigoplus_{k\ge0}C^k_{Z(A)}(\Der(A),A)
$$
is a graded differential subalgebra of $(C(\Der(A),A),d)$, which is
convenient in the induced structure.

Note that in general neither $C(\Der(A),A)$ not $C_{Z(A)}(\Der(A),A)$
are functorial in $A$. Under strong asssumptions on $A$ at least
$C_{Z(A)}(\Der(A),A)$ becomes a
covariant functor in $A$: If $A=C^\infty(M,\Bbb R)$ for
a smooth finite dimensional second countable manifold $M$ we get
$C_{Z(A)}(\Der(A),A)=\Om^*(M)$, the usual graded commutative
differential algebra of usual differential forms.

The convenient structures on $C(\Der(A),A)$ and $C_{Z(A)}(\Der(A),A)$
make them useful as image spaces for constructions.

\subhead\nmb.{5.2}. Insertion operators \endsubhead
For $X\in\Der(A)$ we consider the insertion operator
$$\gather
i_X:C^*(\Der(A),A)\to C^{*-1}(\Der(A),A),\\
(i_X\om)(X_1,\dots,X_k) := \om(X,X_1,\dots,X_k),
\endgather$$
which is a derivation of degree $-1$ satisfying
$i_X|C^0(\Der(A),A)=i_X|A=0$. Note that $i_X$ maps the graded
differential subalgebra $C_{Z(A)}(\Der(A),A)$ into itself.

More generally, let
$K\in C^k(\Der(A),\Der(A)):= L^k_{\text{skew}}(\Der(A);\Der(A))$
be a boun\-ded skew symmetric $k$-linear mapping
$\Der(A)^{k+1}\to \Der(A)$. Then we may define:
$$\align
i_K& : C^l(\Der(A),A) \to C^{l+k-1}(\Der(A),A)\\
(i_K\om)&(X_1,\dots,X_{l+k-1}) =\tag1\\
&= {\tsize\frac1{k!\,(l-1)!}} \sum_{\si\in\Cal S_{k+l-1}}
\sign(\si)\; \om(K(X_{\si1},\dots, X_{\si k}), X_{\si(k+1)},\dotsc)
\endalign$$
Note that $i_{Id_{\Der(A)}}\om = l.\om$ for $\om\in C^l(\Der(A),A)$.

\proclaim{Lemma} $i_K$ is a derivation of degree $k-1$,
$i_K\in\Der_{k-1}(C^*(\Der(A),A))$.
\endproclaim

\demo{Proof}
Clearly one may write
$$\align
&(i_K\om)(X_1,\dots,X_{l+k-1}) = \tag2\\
= \sum_{i_1<\dots<i_k}&
(-1)^{i_1+\dots+i_k-\frac{k(k+1)}{2}} \om(K(X_{i_1},\dots,
X_{i_k}), X_1,\dots,\widehat{X_{i_1}},\dots, \widehat{X_{i_2}},\dotsc).
\endalign$$
Using expression \thetag2 one may check directly that then we have
$$
[i_X,i_K] = i_{K(X,\dots)} = i_{i_XK}\text{ for any }X\in\Der(A).
\tag3$$
 From \thetag3 it is then easy to prove that for
$\ph\in C^p(\Der(A),A)$ and $\ps\in C^q(\Der(A),A)$ we have
$$i_K(\ph.\ps)=(i_K\ph).\ps+(-1)^{(k-1)p}\ph.(i_K\ps),$$
by induction on $k+p+q$.
\qed\enddemo

\subhead\nmb.{5.3}. Remark \endsubhead
In general there exist more bounded graded derivations $D$ in the
space
$\Der(C(\Der(A),A))$ which vanish on $A$, than just those of the form
$i_K$ for some
$K\in C(\Der(A),\Der(A))$.
We shall give an explicit description of all graded derivations on
$C(\Der(A),A)$ in section \nmb!{6} below, but we will not extend the
Fr\"olicher-Nijenhuis calculus to this description: it is too
complicated.

\subhead\nmb.{5.4} \endsubhead
Now let
$K\in C^k(\Der(A),\Der(A))$ and $L\in C^l(\Der(A),\Der(A))$. Then for
the graded commutator we have
$$\align
[i_K,i_L] &= i_K\,i_L -(-1)^{(k-1)(l-1)}i_L\,i_K
     = i_{[K,L]^\wedge}, \text{ where }\\
[K,L]^\wedge &= i_KL - (-1)^{(k-1)(l-1)}i_LK \in C^{k+l-1}(\Der(A),\Der(A))
\endalign$$
is the \idx{\it Nijenhuis-Richardson bracket}, see \cit!{20}, a graded
Lie bracket on the graded vector space $C^{*+1}(\Der(A),\Der(A))$.
Here $i_L$ acts on $C(\Der(A),\Der(A))$ by the same formula as in
\nmb!{5.2}.\thetag1.

\subhead\nmb.{5.5}. Lie derivations \endsubhead
For
$K\in C^k(\Der(A),\Der(A))$ let us now define the
\idx{\it Lie derivation} $\L_K$ along $K$ by the graded commutator
$$\L_K = [i_K,d] = i_K\,d -(-1)^{k-1}d\,i_K\in \Der_k(C(\Der(A),A)).$$
Note that $\L_{Id_{\Der(A)}} = d$.

\proclaim{\nmb.{5.6}. Lemma}
The Lie
derivative of $\om\in C^l(\Der(A),A)$ along $K$ is given by
$$\align
&(\L_K\om)(X_1,\dots,X_{k+l}) = \\
&= {\tsize\frac 1{k!\,l!}}\sum_{\si} \sign\si\;
\L_{K(X_{\si1},\dots,X_{\si k})}(\om(X_{\si(k+1)},\dots,X_{\si(k+l)}))\\
& + (-1)^k\biggl({\tsize\frac{1}{k!\,(l-1)!}}\sum_\si \sign\si\;
\om([X_{\si1},K(X_{\si2},\dots,X_{\si(k+1)})], X_{\si(k+2)},\ldots) \\
&\qquad - {\tsize\frac{1}{(k-1)!\,(l-1)!\,2!}}
\sum_\si \sign \si\; \om(K([X_{\si 1},X_{\si 2}],
X_{\si 3},\ldots), X_{\si(k+2)},\ldots)\biggr)
\endalign$$
\endproclaim

\demo{Proof}
This can be shown by a direct computation starting from formula
\nmb!{5.2}.\thetag2.
\qed\enddemo

\proclaim{\nmb.{5.7}. Proposition}
Let $K\in C^k(\Der(A),\Der(A))$ and $L\in C^l(\Der(A),\Der(A))$. Then
for the graded commutator we have
$$
[\L_K,\L_L] = \L_K\,\L_L -(-1)^{(k-1)(l-1)}\L_L\,\L_K
     = \L_{[K,L]},
$$
where $[K,L] \in C^{k+l}(\Der(A),\Der(A))$
is given by the following formula
$$\align
&[K,L](\row X1{k+l}) = \tag1\\
&= {\tsize\frac 1{k!\,l!}}\sum_{\si} \sign\si\;
[K(X_{\si1},\dots,X_{\si k}),L(X_{\si (k+1)},\dots,X_{\si(k+l)})]\\
& + (-1)^k\biggl({\tsize\frac{1}{k!\,(l-1)!}}\sum_\si \sign\si\;
L([X_{\si1},K(X_{\si2},\dots,X_{\si(k+1)})], X_{\si(k+2)},\ldots) \\
&\qquad - {\tsize\frac{1}{(k-1)!\,(l-1)!\,2!}}
\sum_\si \sign \si\; L(K([X_{\si 1},X_{\si 2}],
X_{\si 3},\ldots), X_{\si(k+2)},\ldots)\biggr) \\
&-(-1)^{kl+l}\biggl({\tsize\frac{1}{(k-1)!\,l!}}
	\sum_\si \sign\si\;
K([X_{\si1},L(X_{\si2},\dots,X_{\si(l+1)}),X_{\si(l+2)},\ldots) \\
&\qquad - {\tsize\frac{1}{(k-1)!\,(l-1)!\,2!}}
\sum_\si \sign \si\; K(L([X_{\si 1},X_{\si 2}],
X_{\si 3},\ldots), X_{\si(l+2)},\ldots)\biggr).
\endalign$$
The bracket $[K,L]$ is called the \idx{\it Fr\"olicher-Nijenhuis
bracket}. It is a graded Lie bracket on $C^*(\Der(A),\Der(A))$.
\endproclaim

\demo{Proof}
The Chevalley coboundary operator for the adjoint representation of
the Lie algebra $\Der(A)$ is given by
$$\align
\partial K(X_0,\dots,X_k) &= \sum_{0\le i\le k}(-1)^i
	[X_i,K(X_0,\dots,\widehat{X_i},\dots,X_k)]\\
& \quad+\sum_{0\le i<j\le k}(-1)^{i+j}
 	K([X_i,X_j],X_0,\dots,\widehat{X_i},\dots,\widehat{X_j},\dots,X_k),
\endalign$$
and it is well known that $\partial\partial=0$.
Then the Fr\"olicher-Nijenhuis bracket \thetag1 is given by
$$
[K,L]=[K,L]_{\wedge } + (-1)^k i(\partial K)L
	-(-1)^{kl+l}i(\partial L)K,\tag2$$
where we have put
$$\multline
[K,L]_\wedge (X_1,\dots,X_{k+l}) = \\
= \tfrac1{k!\,l!}\sum_{\si}\sign(\si)[K(X_{\si1},\dots,X_{\si k}),
     L(X_{\si(k+1)},\dots,X_{\si(k+l)})]_{\Der(A)}.
\endmultline\tag3$$
Formula \thetag2 is the same as in \cit!{20}, p\. 100, where it is also
stated that from this formula `one can show (with a good deal of
effort) that this bracket defines a graded Lie algebra structure'.
Similarly we can write the Lie derivative \nmb!{5.6} as
$$\L_K = \L_\wedge(K) + (-1)^ki(\partial K),\tag4$$
where the action $\L$ of $\Der(A)$ on $A$ is extended to
$\L_\wedge :C(\Der(A),\Der(A))\x C(\Der(A),A)\to C(\Der(A),A)$ by
$$\multline
(\L_\wedge(K)\om)(X_1,\dots,X_{q+k}) =\\
= \tfrac1{k!\,l!}\sum_{\si}\sign(\si)\L(K(X_{\si1},\dots,X_{\si k}))(
     \om(X_{\si(k+1)},\dots,X_{\si(k+q)})).
\endmultline\tag5$$
Using \thetag4 we see that
$$\align
[\L_K,\L_L] &= \L_\wedge(K)\L_\wedge(L)
     -(-1)^{kl}\L_\wedge(L)\L_\wedge(K)\tag6\\
&\quad +(-1)^ki(\partial K)\L_\wedge(L)
     -(-1)^{kl+k}\L_\wedge(L)i(\partial K)\\
&\quad -(-1)^{kl+l}i(\partial L)\L_\wedge(K)
     +(-1)^l\L_\wedge(K)i(\partial L)\\
&\quad +(-1)^{k+l}i(\partial K)i(\partial L)
     -(-1)^{kl+k+l}i(\partial L)i(\partial K),
\endalign$$
and from \thetag2 and \thetag4 we get
$$\align
\L_{[K,L]} &= \L_{[K,L]_\wedge} +(-1)^k\L_{i(\partial K)L}
     -(-1)^{kl+l}\L_{i(\partial L)K} \tag7\\
&= \L_\wedge ([K,L]_\wedge) +(-1)^{k+l}i(\partial[K,L]_\wedge) \\
&\quad +(-1)^k\L_\wedge (i(\partial K)L)
     +(-1)^ki(\partial i(\partial K)L) \\
&\quad -(-1)^{kl+l}\L_\wedge (i(\partial L)K)
     -(-1)^{kl+k}i(\partial i(\partial L)K).
\endalign$$
By a straightforward direct computation one checks that
$$
\L_\wedge(K)\L_\wedge(L)-(-1)^{kl}\L_\wedge(L)\L_\wedge(K) =
\L_\wedge([K,L]_\wedge).
\tag8$$
The expression $\L_\wedge(L)\om$ looks like the `wedge' product in
\nmb!{5.1} to the derivation $i_K$ of degree $k$.
So we may apply lemma \nmb!{5.2} and get
$$
i_K\L_\wedge(L)\om = \L_\wedge(i_KL)\om
     +(-1)^{(k-1)l}\L_\wedge(L)i_K\om
\tag9$$
By a straightforward combinatorial computation one can check directly
from the definitions that the following formula holds:
$$
\partial(i_KL) = i_{\partial K}L +(-1)^{k-1}i_K\partial L
+(-1)^k[K,L]_\wedge \tag{10}
$$
Moreover it is obvious that
$$
\partial[K,L]_\wedge = [\partial K,L]_\wedge +(-1)^k[K,\partial L]_\wedge
\tag{11}$$
We have to show that \thetag6 equals \thetag7. This follows now by
using \thetag8, twice \thetag9, twice \thetag{10}, \thetag{11}, and
$\partial\partial=0$.

That the Fr\"olicher-Nijenhuis bracket defines a graded Lie bracket
follows now from the fact that
$\L:C(\Der(A),\Der(A)) \to \Der(C(\Der(A),A))$ is injective.
\qed\enddemo

\proclaim{\nmb.{5.8}. Lemma} For $K \in C^k(\Der(A),\Der(A))$ and $L \in
C^l(\Der(A),\Der(A))$ we have
$$\align [\L_K,i_L] &= i([K,L]) - (-1)^{k(l-1)}\L(i_LK)\text{, or}\\
[i_L,\L_K] &= \L(i_LK) + (-1)^k\,i([L,K]).\endalign$$
\endproclaim

\demo{Proof}
The two equations are obviously equivalent by graded skew symmetry,
and the second one follows by inserting the following formulas, all from
\nmb!{5.7}: expand the equation by \thetag4, \thetag2, and use then
\thetag{10}.
\qed\enddemo

\subhead\nmb.{5.9}. Remark \endsubhead
As formal consequences of lemma \nmb!{5.8} we get
the following fomulae:
For
$K_i \in C^{k_i}(\Der(A),\Der(A))$ and $L_i \in C^{k_i+1}(\Der(A),\Der(A))$ we
have
$$\align
&\left[\L_{K_1}+i_{L_1},\L_{K_2}+i_{L_2}\right] = \tag1 \\
&\qquad = \L\left([K_1,K_2]
    + i_{L_1}K_2 - (-1)^{k_1k_2}i_{L_2}K_1\right)\\
&\qquad\quad + i\left([L_1,L_2]^{\wedge}
    + [K_1,L_2] -(-1)^{k_1k_2}[K_2,L_1]\right).
\endalign$$
Each summand of this formula looks like a semidirect product of
graded Lie algebras, but the mappings
$$\align
i:C(\Der(A),\Der(A)) &\to \End(C(\Der(A),\Der(A)),[\quad,\quad]) \\
\ad:C(\Der(A),\Der(A)) &\to
     \End(C(\Der(A),\Der(A)),[\quad,\quad]^{\wedge})
\endalign$$
do not take values in the subspaces of graded derivations. We
have instead for $K \in C^k(\Der(A),\Der(A))$ and $L \in
C^{l+1}(\Der(A),\Der(A))$ the following relations:
$$\align
&i_L[K_1,K_2] = [i_LK_1,K_2]
     + (-1)^{k_1l}[K_1,i_LK_2] \tag2 \\
&\qquad -\left((-1)^{k_1l}i([K_1,L])K_2
      - (-1)^{(k_1+l)k_2}i([K_2,L])K_1\right)\\ \allowdisplaybreak
&[K,[L_1,L_2]^{\wedge}] = [[K,L_1],L_2]^{\wedge}
      + (-1)^{kk_1}[L_1,[K,L_2]]^{\wedge} -\tag3 \\
&\qquad -\left((-1)^{kk_1}[i(L_1)K,L_2]
      - (-1)^{(k+k_1)k_2}[i(L_2)K,L_1]\right)
\endalign$$
The algebraic meaning of these relations and its
consequences in group theory have been investigated in
\cit!{19}. The corresponding product of groups is well
known to algebraists under the name `Zappa-Szep'-product.

Moreover the Chevalley coboundary operator is a homomorphism from the
Fr\"o\-li\-cher\--Ni\-jen\-huis bracket to the Nijenhuis-Richardson bracket:
$$
\partial[K,L] = [\partial K,\partial L]^\wedge
\tag4$$

\demo{Proof}
\thetag1 follows from repeated applications of \nmb!{5.8}. It is easy
to show that the images under  $\L$ of both sides of \thetag2
coincide, using \nmb!{5.8}, then we jus note that $\L$ is injective.
For \thetag3 it is again easy to show that the images under $i$ of
both sides coincide, again using \nmb!{5.8}; also
$i:C(\Der(A),\Der(A))\to \Der(C(\Der(A),A))$ is
injective.

\thetag4 follows from \nmb!{5.7},\thetag2, \thetag{10}, and
\thetag{11}, and from \nmb!{5.4}.
\qed\enddemo

\subhead\nmb.{5.10}. The graded differential subalgebra
$C_{Z(A)}(\Der(A),\Der(A))$ \endsubhead
Clearly a derivation $i_K\in\Der(C(\Der(A),A))$ for
$K\in C^k(\Der(A),\Der(A))$ maps the graded differential subalgebra
$C_{Z(A)}(\Der(A),A)$ into itself if and only if $K$ is skew
$Z(A)$-multilinear, i\.e\.
$$\gather
K\in C_{Z(A)}(\Der(A),\Der(A))\\
K(z_1X_1,\dots,z_kX_k)=z_1\dots z_kK(X_1,\dots,X_k)
     \text{ for all }z_i\in Z(A)\text{ and }X_i\in \Der(A).
\endgather$$
Then also the Lie derivation $\L_K$ respects $C_{Z(A)}(\Der(A),A)$,
and so the closed linear subspace $C_{Z(A)}(\Der(A),\Der(A))$ is a
graded Lie subalgebra for both brackets $[\quad,\quad]^\wedge $ and
$[\quad,\quad]$, and all formulas in this section continue to hold.
But note that in the simpler formula \nmb!{5.7}.\thetag2, although
$[K,L]$ is $Z(A)$-multilinear, none of the three summands is in
$C_{Z(A)}(\Der(A),\Der(A))$. The same applies to
\nmb!{5.7}.\thetag4.

\subhead\nmb.{5.11}. Remark \endsubhead
If $K\in C^1(\Der(A),\Der(A)) = L(\Der(A),\Der(A))$ then formula
\nmb!{5.7}.\thetag1 boils down to
$$
[K,K](X,Y) = 2([KX,KY]-K[KX,Y]-K[X,KY]+K^2[X,Y]).
$$
So if $[K,K]=0$ then the image of $K$ is a
Lie subalgebra of $\Der(A)$. If $K$ is moreover a projection,
$K\o K=K$, then also the kernel is a Lie subalgebra.
In this case one could view $K$ as a `connection' and could say
that the kernel and the
image of $K$ are `involutive subbundles' whose `curvatures
vanish', compare with \cit!{3}, section 5. We shall elaborate on this
topic in a later paper.

If $K:\Der(A)\to \Der(A)$ is semisimple with $[K,K]=0$, then
each eigenspace of $K$ is also a Lie subalgebra. For if $KX=\la X$
and $KY=\la Y$ then we get $(K-\la)^2[X,Y]=0$ and thus
$K[X,Y]=\la[X,Y]$. In particular the kernel of $K$ is a Lie
subalgebra $\la=0$.

If moreover $K$ is $Z(A)$-linear,
$K\in C^1_{Z(A)}(Der(A), Der(A))=Der(A,\underline{\Omega}^1_{Der}(A))$,
then the above Lie subalgebras of $Der(A)$ are also
$Z(A)$-submodules.

\head\totoc\nmb0{6}. Description of all derivations 
in the Chevalley differential graded algebra \endhead

\subhead\nmb.{6.1}. Insertion operators \endsubhead
The space $\Hom^A_A(C^p(\Der(A),A),A)$ is a sort of `$A$-dual' of
$C^p(\Der(A),A)$.
Let us write $\langle \Xi,\ph\rangle_A\in A$ for the evaluation of
the element
$\Xi\in \Hom^A_A(C^p(\Der(A),A),A)$ on $\ph\in C^p(\Der(A),A)$. Note
that for $a,b\in A$ we have
$\langle \Xi,a.\ph.b\rangle_A= a.\langle \Xi,\ph\rangle_A.b$.
Then we consider the (closure of the) linear subspace
$$
\sum_{0<i<p}C^i(\Der(A),A).C^{p-i}(\Der(A),A)\subset C^p(\Der(A),A)
$$
and its annihilator
$$\align
\operatorname{Ann}^p(\Der(A),A):&=
     \left( \sum_{0<i<p}C^i(\Der(A),A).C^{p-i}(\Der(A),A)\right)^{\bigcirc}\\
&\subseteq \Hom^A_A(C^p(\Der(A),A),A),
\endalign$$
and set
$$\align
\operatorname{Ann}^1(\Der(A),A)&= \Hom^A_A(C^1(\Der(A),A),A),\\
\operatorname{Ann}^0(\Der(A),A)&= \Der(A).\\
\endalign$$
For $\Xi\in \operatorname{Ann}^p(\Der(A),A)$ we consider the `insertion
operator'
$$\gather
i_\Xi:C^*(\Der(A),A)\to C^{*-p}(\Der(A),A),\\
(i_\Xi\om)(X_1,\dots,X_l) := \langle\Xi,\om(
     \undersetbrace{p\text{ times}}\to{\quad,\dots,\quad},
     X_1,\dots,X_l)\rangle_A,
\endgather$$
which is a derivation of degree $-p$ satisfying
$i_\Xi|C^q(\Der(A),A)=0$ for $q<p$, since $\Xi$ is an $A$-bimodule
homomorphism and annihilates all `small' products.

More generally let
$K\in C^k(\Der(A),\operatorname{Ann}^p(\Der(A),A))$
be a bounded skew symmetric $k$-linear mapping
$\Der(A)^{k}\to \operatorname{Ann}^p(\Der(A),A)$. Then we define the bounded
linear mapping
$i_K : C^l(\Der(A),A) \to C^{l+k-p}(\Der(A),A)$ by
$i_K |C^q(\Der(A),A)=0$ for $q<p$, and by
$$\multline
(i_K\om)(X_1,\dots,X_{l+k-p}) = \\
= {\tsize\frac1{k!\,(l-1)!}} \sum_{\si\in\Cal S_{k+l-1}}
     \sign\si. \langle K(X_{\si1},\dots, X_{\si k}),
     \om(\undersetbrace{p\text{ times}}\to{\quad,\dots,\quad},
     X_{\si(k+1)},\dotsc)\rangle_A.
\endmultline\tag1$$

\proclaim{Lemma} $i_K$ is a derivation of degree $k-p$,
$i_K\in\Der_{k-p}(C^*(\Der(A),A))$.
\endproclaim

\demo{Proof}
Clearly one may write
$$\multline
(i_K\om)(X_1,\dots,X_{l+k-1}) =\\
     = \sum_{i_1<\dots<i_k} (-1)^{i_1+\dots+i_k-\frac{k(k+1)}{2}}
     \langle K(X_{i_1},\dots,X_{i_k}),\\
\om(\undersetbrace{p\text{ times}}\to{\quad,\dots,\quad},
     X_1,\dots,\widehat{X_{i_1}},\dots,\widehat{X_{i_2}},\dotsc)\rangle_A.
\endmultline\tag2$$
Using expression \thetag2 one may check directly that then we have
$$
[i_X,i_K] = i_{K(X,\dots)} = i_{i_XK}\text{ for any }X\in\Der(A).
\tag3$$
 From \thetag3 it is easy to prove that for
$\ph\in C^m(\Der(A),A)$ and $\ps\in C^n(\Der(A),A)$ we have
$$i_K(\ph.\ps)=(i_K\ph).\ps+(-1)^{(k-p)m}\ph.(i_K\ps),$$
by induction on $k-p+m+n$, for each fixed $p$.
\qed\enddemo

\proclaim{\nmb.{6.2}. Proposition}
Let $D\in\Der_k(C(\Der(A),A))$ be a graded derivation.
Then there are unique
$K_0\in C^k(\Der(A),\Der(A))=C^k(\Der(A),\operatorname{Ann}^0(\Der(A),A))$ and
$K_p\in C^{k+p}(\Der(A),\operatorname{Ann}^p(\Der(A),A))$ for $p=1,2,\dotsc$
such that
$$D= \L_{K_0} + \sum_{p=1}^\infty i_{K_p}.$$
\endproclaim

Note that on each fixed component $C^l(\Der(A),A)$ this is a finite
sum.

\demo{Proof}
Let $D\in\Der_k(C(\Der(A),A))$ be a graded derivation. Then the
restriction $D|A$ of $D$ to $A=C^1(\Der(A),A)$ is an element of
$$\align
K_0 := D|A &\in \Der(A,C^k(\Der(A),A)) = \\
&= \Der(A,L^k_{\text{skew}}(\Der(A),A)) =\\
&= L^k_{\text{skew}}(\Der(A),\Der(A,A)),\text{ by \nmb!{2.13}} \\
&= C^k(\Der(A),\Der(A)).
\endalign$$
By \nmb!{5.5} we have the Lie derivation
$\L_{K_0}\in \Der_k(C(\Der(A),A))$ which coincides with $D$ on
$A=C^0(\Der(A),A)$, so that the difference
$D-\L_{K_0}\in\Der_k(C(\Der(A),A))$ is a graded derivation with
$(D-\L_{K_0})|A=0$.
So
$$\align
K_1:=(D-\L_{K_0})|C^1(\Der(A), &A)\in
\Hom^A_A(C^1(\Der(A),A),C^{k+1}(\Der(A),A)) \\
&= \Hom^A_A(L(\Der(A),A),L^{k+1}_{\text{skew}}(\Der(A),A)) \\
&= L^{k+1}_{\text{skew}}(\Der(A),\Hom^A_A(L(\Der(A),A),A))
     \text{ by \nmb!{2.13}} \\
&= L^{k+1}_{\text{skew}}(\Der(A),\operatorname{Ann}^1(\Der(A),A))\\
\endalign$$
By \nmb!{6.1} we get a derivation $i_{K_1}:C(\Der(A),A)
\to C(\Der(A),A)$ and the difference
$D-\L_{K_0}-i_{K_1}\in\Der_k(C(\Der(A),A))$ now vanishes on
$A=C^0(\Der(A),A)$
and $C^1(\Der(A),A)$. Thus the restriction is an $A$-bimodule
homomorphism
$$
(D-\L_{K_0}-i_{K_1})|C^2(\Der(A),A):C^2(\Der(A),A)\to C^{k+2}(\Der(A),A))
$$
vanishes on all products of 1-forms. If we consider
$$\align
K_2:=(D-\L_{K_0}-i_{K_1})&|C^2(\Der(A),A)\in
     \Hom^A_A(C^2(\Der(A),A),C^{k+2}(\Der(A),A)) \\
&= \Hom^A_A(C^2(\Der(A),A),L^{k+2}_{\text{skew}}(\Der(A),A)) \\
&= L^{k+2}_{\text{skew}}(\Der(A),\Hom^A_A(C^2(\Der(A),A),A))\text{ by
     \nmb!{2.13}},\\
&= C^{k+2}(\Der(A),\operatorname{Ann}^2(\Der(A),A)),
\endalign$$
then by lemma \nmb!{6.1} we have a derivation
$i_{K_2}\in \Der_k(C(\Der(A),A))$ which vanishes on $C^q(\Der(A),A)$
for $q=0,1$ and coincides with $D-\L_{K_0}-i_{K_1}$ on
$C^2(\Der(A),A)$, so the derivation
$D-\L_{K_0}-i_{K_1}-i_{K_2}$ vanishes on $C^q(\Der(A),A)$ for $q=0,1,2$, and
we may repeat the process.
\qed\enddemo

\subhead\nmb.{6.3}. Remarks \endsubhead
If we try to expand the graded commutator $[i_K,i_L]$ using
\nmb!{6.2}, the resulting formulas are too complicated to be written
down easily: we should invent new notation. We refrain from doing
this since we also have no use for it.

For $K\in C^k(\Der(A),\operatorname{Ann}^1(\Der(A),A))$ we get
$[i_K,d]|A = \L_{\hat K}|A$, where $\hat K\in C^k(\Der(A),\Der(A))$
is given by
$\hat K(X_1,\dots,X_k)(a) = K(X_1,\dots,X_k)(da)$.
The higher order part in the expansion of $[i_K,d]$ according to
\nmb!{6.2} does not vanish, and it can be written down in principle.

If we want to classify all derivations in $\Der(C_{Z(A)}(\Der(A),A))$
we can just repeat the developpment in this section, but have to  replace
$\operatorname{Ann}^k(\Der(A),A)$ everywhere by the annihilator
$$\align
\operatorname{Ann}^p_{Z(A)}(\Der(A),A):&=
     \left( \sum_{0<i<p}C^i_{Z(A)}(\Der(A),A).
     C^{p-i}_{Z(A)}(\Der(A),A)\right)^{\bigcirc}\\
&\subseteq \Hom^A_A(C^p_{Z(A)}(\Der(A),A),A),\\
\operatorname{Ann}^1_{Z(A)}(\Der(A),A)&= \Hom^A_A(C^1_{Z(A)}(\Der(A),A)),\\
\operatorname{Ann}^0_{Z(A)}(\Der(A),A)&= \operatorname{Ann}^0(\Der(A),A)
=\Der(A).\\
\endalign$$

\head\totoc\nmb0{7}. Diagonal bimodules \endhead

\subhead\nmb.{7.1}. The differential graded algebra
$\Om_\Der(A)$ \endsubhead
By the universal property of $(\Om^*(A),d)$ there is a unique
bounded homomorphism of graded differential algebras
$$\CD
A @>d>> \Om^1(A) @>d>> \Om^2(A) @>d>> \cdots\\
@|               @VV\ze_1V      @VV{\ze_2}V    @.\\
A @>d>> C^1_{Z(A)}(\Der(A),A) @>d>> C^2_{Z(A)}(\Der(A),A) @>d>> \cdots
\endCD$$
which is given by
$$(\ze_k\om)(X_1,\dots,X_k) = j_{X_k}\dots j_{X_1}\om \in \Om^0(A)=A,$$
for $\om\in\Om^k(A)$ and $X_i\in\Der(A)$.
The kernel of $\ze$ is the space
$$\align
F^1\Om^*(A) &= \bigoplus_{k\ge0}F^1\Om^k(A)\\
F^1\Om^k(A) &= \{\om\in\Om^k(A): j_{X_1}\dots j_{X_k}\om=0\text{ for
all }X_i\in\Der(A)\},
\endalign$$
see \cit!{6}, which is a closed graded differential ideal in
$(\Om(A),d)$ by a short computation. It is part of an obvious
filtration which leads to a spectral sequence, see also \cit!{6}.
The image of the homomorphism $\ze$ is denoted by
$(\Om_\Der^*(A),d)$
and it will be equipped with with the quotient structure of a
convenient vector space, which is a finer structure than that induced
from $C_{Z(A)}(\Der(A),A)$.

Note that $\Om_{\text{Der}}(A)$ is not functorial in $A$ in general;
its convenient structure makes it useful as a source for
constructions.

\subhead\nmb.{7.2}. Derivation based bimodules \endsubhead
The space $\Om^1_\Der(A) =
\Om^1(A)/\{\om\in \Om^1(A):j_X\om=0\text{ for all }X\in\Der(A)\}$
has the following remainder of the unversal property
\nmb!{3.3}.\therosteritem1:
\roster
\item For any bounded derivation $X:A\to A$
     there is a unique bounded $A$-bimodule homomorphism
     $i_X:\Om^1_\Der(A)\to A$
     such that $X=i_X\o d$.
\endroster
We now say that an $A$-bimodule $M$ is a \idx{\it derivation-based
$A$-bimodule} if the following property is satisfied:
\roster
\item[2] For any bounded derivation $D:A\to M$
     there is a bounded $A$-bimodule homomorphism
     $i_D:\Om^1_\Der(A)\to M$ such that $D=i_D\o d$. In fact, if
     $i_D$ exists it is unique since the image of $d$ generates
     $\Om^1_{\text{Der}}(A)$ as $A$-bimodule.
\endroster
By the universal property \nmb!{3.3}.\therosteritem1 of $\Om^1(A)$
condition \therosteritem2 is equivalent to the following:
\roster
\item [3] $\Hom^A_A(\Om^1(A),M) = \Hom^A_A(\Om^1_\Der(A),M)$.
\endroster
It is obvious that for an arbitrary index set $J$ the direct product
$A^J$ with the product $A$-bimodule structure has property
\therosteritem2, and also each of it's sub $A$-bimodules.
Let us call a \idx{\it diagonal bimodules}
any $A$-bimodule which is isomorphic to a submodule of some product
$A^J$. So all diagonal bimodules are derivation based.

\proclaim{\nmb.{7.3}. Proposition}
An $A$-bimodule $K$ is diagonal if and only if the `$A$-dual'
$\Hom^A_A(K,A)$ separates points on $K$. For each $A$-bimodule $M$
there exists a universal diagonal quotient $p_M:M\to \Diag(M)$ such
that each bimodule homomorphism from $M$ into a diagonal bimodule $K$
factors over $\Diag(M)$:
$$\CD
M @>>> K \\
@V{p_M}VV @|\\
\Diag(M) @>>> K  &
\endCD$$
Let $\underline{\operatorname{Bimod}}$ be the category of all $A$-bimodules
and let $\underline{\Diag}$ denote the full subcategory of
all diagonal $A$-bimodules, with
$\io:\underline{\Diag}\to \underline{\operatorname{Bimod}}$ the
embedding and
$\Diag:\underline{\operatorname{Bimod}}:\to \underline{\Diag}$ the
functor from above. Then $\Diag$ is left adjoint to $\io$, i\.e\. we
have the following natural correspondence:
$$\Hom^A_A(\Diag(M),K) \cong \Hom^A_A(M,\io K)$$
Thus the functor $\Diag$ respects colimits, whereas $\io$ respects
limits, and the category $\underline{\Diag}$ is complete: products
and submodules of modules in $\underline{\Diag}$ are again in
$\underline{\Diag}$.
\endproclaim

\demo{Proof}
The vector space $\Hom^A_A(A^J,A)$ clearly separates points on $A^J$,
so this is also true for each sub bimodule of $A^J$.
For an arbitrary  $A$-bimodule $M$ we consider the following
homomorphism of $A$-bimodules:
$$
M @>p_M>> A^{\Hom^A_A(M,A)}, \quad
M\ni m \mapsto (\ph(m))_{\ph\in \Hom^A_A(M,A)} \in A^{\Hom^A_A(M,A)},
$$
and we denote by $\Diag(M)$ the image of $p_M$, an diagonal $A$-bimodule.
If the vector space $\Hom^A_A(M,A)$ separates points on $M$ then $p_M$ is
injective and
$M$ is diagonal. The kernel of $p$ is given by
$\ker(p_M)=\{m\in M: \ph(m)=0\text{ for all }\ph\in \Hom^A_A(M,A)\}$,
and obviously any homomorphism $M\to K$ into a diagonal $A$-bimodule
$K$ vanishes on $\ker(p_M)$ and thus factors over $p_M$.

The remaining statements follow by basic category theory.
\qed\enddemo

\subhead\nmb.{7.4}. Examples \endsubhead
We have $\Diag(\Om^1(A)) = \Om^1_\Der(A)$, since any homomorphism
$\Om^1(A)\to A$ corresponds to a derivation $A\to A$ and thus factors
to a homomorphism $\Om^1_\Der(A)\to A$.

Obviously we have $\Diag(A)=A$, but let us consider $\Diag(A\tilde\otimes A)$.
We have $\Hom^A_A(A\tilde\otimes A,A)\cong A$ by
$\ph\mapsto \ph(1\otimes 1)$, so
$\ker(p_{A\tilde\otimes A}:A\tilde\otimes A\to\Diag(A\tilde\otimes A))$
consists of all
\roster
\item $\sum_n a_n\tilde\otimes b_n \in A\tilde\otimes A$ which
       satisfy $\sum_n a_n.c.b_n=0$ for each $c\in A$. \endroster
Taking $c=1$ in \therosteritem1 we see that
$\ker(p_{A\otimes A})\subset \Om^1(A)$, and for each element in
\therosteritem1 we have
$$\gather
\sum_n a_n\tilde\otimes b_n =\sum_n a_n\,db_n,\\
i_{\ad(c)}\sum_n a_n\,db_n = \sum_n a_n.(c.b_n-b_n.c)=0,
\endgather$$
so that $\ker(p_{A\otimes A})$ is the subspace $\Om^1(A)_{H-\Int(A)}$
of all elements in $\Om^1(A)$ which are horizontal with respect to
$\Int(A)$, see \nmb!{3.5}.

We also see that $1\otimes a - a\otimes 1 \in \ker(p_{A\tilde\otimes A})$
if and only if $a\in Z(A)$, so that $p_{A\tilde\otimes A}$ factors as
follows:
$$
A\tilde\otimes A \to A\tilde\otimes_{Z(A)} A
     \to \Diag(A\tilde\otimes A).
\tag2$$
Since $\Diag$ is a left adjoint functor, it is right exact, and we
get the following diagram with exact rows and exact columns:
$$\CD
@. 0 @. 0 @. @. \\
@. @VVV @VVV @. @. \\
0 @>>> F^1\Om^1(A) @>>> \Om^1(A)_{H-\Int(A)} @>>>
     \Om^1_\Der(A)_{H-\Int(A)} @>>> 0\\
@. @VVV @VVV @V0VV @.\\
0 @>>> \Om^1(A) @>>> A\tilde\otimes A @>>> A @>>> 0 \\
@. @V{p_{\Om^1(A)}}VV @V{p_{A\tilde\otimes A}}VV @| @.\\
@. \Om^1_{\Der}(A) @>>> \Diag(A\tilde\otimes A) @>>> A @>>> 0 \\
@. @VVV @VVV @VVV @.\\
@. 0 @. 0 @. 0 @.\\
\endCD$$
 From this and the factorization \thetag2 it follows, that for
commutative $A$ we have $\Diag(A\tilde\otimes A)= A$.

\proclaim{\nmb.{7.5}. Lemma} Over the algebra
$\operatorname{Mat}_N(\Bbb C)$ of complex $(N\x N)$-matrices every
bimodule is diagonal.
\endproclaim

\demo{Proof}
It is well known that every irreducible left
$\operatorname{Mat}_N$-module is isomorphic to $\Bbb C^N$, and that
every left $\operatorname{Mat}_N$-module is semisimple. From this by
transfinite induction one can show that every left
$\operatorname{Mat}_N$-module $L$ (with its finest locally convex
topology) is isomorphic to the direct sum of copies of $\Bbb C^N$,
$L\cong \Bbb C^N\otimes E$ for a vector space $E$.

Now let $M$ be a $\operatorname{Mat}_N$-bimodule. As a left
$\operatorname{Mat}_N$-module we have $M\cong \Bbb C^N\otimes E$. Take a
minimal projection $p\in \operatorname{Mat}_N$, then
$p.M\cong \Bbb C.v\otimes E\cong E$ is a right
$\operatorname{Mat}_N$-module. By the argument above, as a right
$\operatorname{Mat}_N$-module we have
$E\cong K\otimes (\Bbb C^N)^*$. Thus
$$
M\cong C^N\otimes K\otimes (\Bbb C^N)^*
\cong \Bbb C^N\otimes (\Bbb C^N)^*\otimes K
\cong \operatorname{Mat}_N\otimes K
\subset \operatorname{Mat}_N{}^K.\qed
$$
\enddemo

\head\totoc\nmb0{8}. Derivations on the differential graded algebra
$\Om_{\text{Der}}(A)$ \endhead

\proclaim{\nmb.{8.1}. Theorem}
Let $A$ be a convenient algebra. Then we have the following bounded
canonical mappings, where on the left hand side all mappings are
homomorphisms of graded differential algebras
$$\minCDarrowwidth{15pt}\CD
\Om^*(A)  @. \qquad\qquad @.           \Der(A,\Om^*(A)) @. \\
@VV{\ze}V     @.          @VV{\ze}V   @.  \\
\Om^*_\Der(A) @. \qquad\qquad @.          \Der(A,\Om^*_\Der(A)) @.\\
@V{\subset}V{\bar\ze}V    @.           @V{\subset}V{\bar\ze}V  @.   \\
C_{Z(A)}(\Der(A),A) @. \qquad\qquad @.    \Der(A,C_{Z(A)}(\Der(A),A))
                              @= C_{Z(A)}(\Der(A),\Der(A))\\
@V{\subset}VV      @.         @V{\subset}VV  @.   \\
C(\Der(A),A) @. \qquad\qquad @.    \Der(A,C(\Der(A),A))
                              @= C(\Der(A),\Der(A))
\endCD$$
Then for every element $K$ of degree $k$ in one of the right hand
spaces there is a canonical graded derivation $i_K$ (resp\. $j_K$)
of degree $k-1$ on the
corresponding left hand space which vanishes in degree 0;
this corresponds to the graded Lie bracket $[\quad,\quad]^\wedge $ on
the right hand spaces. There is also the corresponding Lie derivation
$\L_K$ of degree $k$ on the left hand space, which leads to the
Fr\"olicher-Nijenhuis bracket $[\quad,\quad]$ on the right hand space.
The vertical arrows intertwine all this derivations and the right
hand ones are homomorphisms for all brackets mentioned.
\endproclaim

The proof of this theorem will fill the rest of this section
\nmb!{8}.

\subhead\nmb.{8.2} \endsubhead
Let $D\in \Der_k(\Om^*_\Der(A))$ be a derivation which
vanishes on $A=\Om^0_\Der(A)$. Then
by setting $K := (D|\Om^1_\Der(A))\o d \in \Der(A,\Om^{k+1}_\Der(A))$
we have
$$D(\om_0\dots\om_l) =
\sum (-1)^{ik}(\om_0)\dots (\om_i)\dots(\om_l),\tag1$$
so $D$ is uniquely determined by $D|\Om^1_\Der(A)$, so by $K$.

Conversely, let $K\in\Der(A,\Om^{k+1}_\Der(A))$ and consider the
corresponding homomorphism
$\tilde J_K\in \Hom^A_A(\Om^1(A),\Om^{k+1}_\Der(A))$ with
$\tilde j_K\o d = K$.
We extend it to $\tilde j_K:\Om^*(A)\to \Om^{*+k}_\Der(A)$
by the right hand side of the universal analogon of \thetag1, i\.e\.
$$\tilde j_K(\om_0\otimes_A\dots\otimes\om_l) =
\sum (-1)^ik\ze(\om_0)\dots \tilde j_K(\om_i)\dots\ze(\om_l)\tag2$$
Then $\tilde j_K$ is a graded derivation of degree $k$ along $\ze$.
We are going to show that $\tilde j_K$ factors to
$\Om^*_\Der(A)$, but we need some preparation.

\subhead\nmb.{8.3} \endsubhead
For $X\in \Der(A) \cong \Hom^A_A(\Om^1(A),A)
\cong \Hom^A_A(\Om^1_\Der(A),A)$ we clearly have the
following factorization:
$$\CD
\Om^*(A) @>j_X>> \Om^{*-1}(A)\\
@V\ze VV          @V\ze VV\\
\Om^*_\Der(A) @>i_X>> \Om^{*-1}_\Der(A)\\
@VVV          @VVV\\
C^*_{Z(A)}(\Der(A),A) @>i_X>> C^*_{Z(A)}(\Der(A),A)
\endCD$$
For $K\in \Der(A,\Om^{k+1}_\Der(A))$ let us
consider the `graded commutator'
$$i_X\,\tilde j_K - (-1)^k\tilde j_K\,j_X: \Om^*(A)
     \to \Om^{*+k}_\Der(A)$$
It is still a graded derivation along $\ze$, and by applying
\nmb!{8.2} we see that
$$\gather
i_X\,\tilde j_K - (-1)^k\tilde j_K\,j_X = \tilde j_{[X,K]^\De}\\
\text{ for }[X,K]^\De\in \Der(A,\Om^{k+1}_\Der(A)),\\
[X,K]^\De(\om)= i_XK(\om)-(-1)^k\tilde j_K(j_X\om)=(i_X\o K)(\om).
\endgather$$

\proclaim{\nmb.{8.4}. Lemma} For
$K\in \Der(A,\Om^{k+1}_\Der(A))$ we have
$\tilde j_K(F^1(\Om^l(A)))=0$ for all $l$.
\endproclaim

\demo{Proof}
We do induction on $l+k$. For $l=0$ we have $F^1(A)=0$, so the
assertion holds. If $k=-1$ then $K=X\in \Der(A)$ and we have
$\tilde j_X(F^1(\Om^l(A))) =
\ze(j_X(F^1(\Om^l(A))))\subset \ze(F^1(\Om^l(A))) =0$.

Now for the induction step we take $\ph\in F^1(\Om^l(A))$, so that
$j_{X_1}\dots j_{X_l}\ph=0$ for all $X_i\in \Der(A)$. From \nmb!{8.3}
we get
$$
i_X\,\tilde j_K\,\ph = (-1)^k \tilde j_K(j_X\ph) +
\tilde j_{i_X\o K}\ph = 0,
$$
by induction:
the first summand vanishes since $j_X\ph\in F^1\Om^{l-1}(A)$, and the
second summand vanishes since
$i_X\o K\in \Der(A,\Om^{k}_\Der(A))$.
\qed\enddemo

\proclaim{\nmb.{8.5}. Corollary} For each $k$ the $A$-bimodule
$\Om^k_\Der(A)$ is a derivation based bimodule.
\endproclaim

\demo{Proof}
By the unversal property \nmb!{3.3}.\therosteritem1 and from
\nmb!{8.4} we get
$$\Der(A,\Om^k_\Der(A))
\cong \Hom^A_A(\Om^1(A),\Om^k_\Der(A)) @<{\ze^*}<\cong<
\Hom^A_A(\Om^1_\Der(A),\Om^k_\Der(A)).\qed$$
\enddemo

This result also follows from the fact that $\Om^k_\Der(A)$ is a sub
$A$-bimodule of a direct product $A^J$, i\.e\. a diagonal bimodule,
see the last remark in
\nmb!{7.2}. By the same reason the bimodules $C^k(\Der(A),A)$ and
$C^k_{Z(A)}(\Der(A),A)$ are also diagonal bimodules, thus
derivation based $A$-bimodules.

\subhead\nmb.{8.6}. Insertion operators \endsubhead
For any $K\in \Der(A,\Om^{k+1}_\Der(A))$ we may now define
the insertion operator $i_K\in \Der(\Om^*_\Der(A))$ by the
the following factorization which is due to lemma \nmb!{8.4}:
$$\CD
\Om^*(A) @>\tilde j_K>> \Om^{*+k}_\Der(A)\\
@V\ze VV              @|\\
\Om^*_\Der(A)  @>i_K>> \Om^{*+k}_\Der(A)\\
\endCD$$
 From \nmb!{8.2} we may now conclude that any derivation
$D\in \Der(\Om^*_\Der(A))$ with $D|A=0$ is of the form $i_K$
for a unique $K\in \Der(A,\Om^{k+1}_\Der(A))$.

For $K,L\in \Der(A,\Om^{*}_\Der(A))$ of degree $k+1$
and $l+1$, respectively, we have
$$\gather
[i_K,i_L]=i_K\,i_L -(-1)^{kl} i_L\,i_K = i_{[K,L]^\wedge },\text{ where}\\
[K,L]^\wedge = i_K\o L -(-1)^{kl} i_L\o K:
     \Om^1(A)\to \Om^{k+l+1}_\Der(A)
\endgather$$

\subhead\nmb.{8.7}. Lie derivations \endsubhead
For $K\in \Der(A,\Om^{k}_\Der(A))$ we define the
\idx{\it Lie derivation} along $K$ by
$$\L_K=[i_K,d]=i_K\,d-(-1)^{k-1}d\,i_K\in\Der_k(\Om^*_\Der(A)).$$
Similar as in \nmb!{4.3} one sees that $\L_{Id}=d$.

\proclaim{\nmb.{8.8}. Proposition}
For any graded derivation $D\in \Der_k(\Om^*_\Der(A))$ there
are unique
$K\in \Der(A,\Om^{k}_\Der(A))$ and
$L\in \Der(A,\Om^{k+1}_\Der(A))$ such
that
$$D=\L_K+i_L.$$
We have $L=0$ if and only if $[D,d]=0$; and $K=0$ if and only if
$D|A=0$.
\endproclaim

\demo{Proof}
$D|A:A=\Om^0_\Der(A)\to \Om^k_\Der(A)$ is a
derivation with values in the derivation based $A$-bimodule
$\Om^k_\Der(A)$ (\nmb!{7.2}), so by the universal property
\nmb!{3.3}.\therosteritem1 and by \nmb!{8.5} there is a
unique
$K\in \Der(A,\Om^{k}_\Der(A))$
with $D|A=K$. But then $(D-\L_K)|A=0$, so by \nmb!{8.6} we have
$D-\L_K=i_L$ for a unique
$L\in \Der(A,\Om^{k+1}_\Der(A))$.
\qed\enddemo

\subhead\nmb.{8.9}. The Fr\"olicher-Nijenhuis bracket \endsubhead
For $K\in \Der(A,\Om^{k}_\Der(A))$ and
$L\in \Der(A,\Om^{l}_\Der(A))$ the
graded commutator of the Lie derivations $[\L_K,\L_L]$ commutes with
$d$, so by \nmb!{8.8} we have
$[\L_K,\L_L]=\L_{[K,L]}$ for unique
$[K,L]\in \Der(A,\Om^{k+l}_\Der(A))$.
We may conclude that this bracket
$[\quad,\quad]$ defines a graded Lie algebra structure on
$\Der(A,\Om^{*}_\Der(A))$, because the mapping
$\L$ is injective.
This bracket is called the \idx{\it Fr\"olicher-Nijenhuis bracket}.

\subhead\nmb.{8.10}  \endsubhead
For $K\in \Der(A,\Om^{k}_\Der(A))$ we define
$\ze K\in C^{k}_{Z(A)}(\Der(A),\Der(A))$ by
$$
\ze K(X_1,\dots,X_k) := i_{X_k}\o\dots\o i_{X_1}\o K:A\to A
$$
If we denote for the moment by
$\bar\ze:\Om^*_\Der(A) \to C^*_{Z(A)}(\Der(A),A)$
the embedding of graded differential algebras, then we have
$$
\bar\ze\o i_K = i_{\ze K}\o\bar\ze:
     \Om^*_\Der(A)\to C^{*+k-1}_{Z(A)}(\Der(A),A),
$$
So the elements
$L\in C^k_{Z(A)}(\Der(A),\Der(A))$ for
which $i_L$ maps the graded subalgebra $\Om^*_\Der(A)$ into itself
are precisely those of the form $L=\ze K$ for some
$K\in \Der(A,\Om^{k}_\Der(A))$.

\proclaim{\nmb.{8.11}. Proposition} The injective bounded linear mapping
$$\ze{}: \Der(A,\Om^{*}_\Der(A))
\to C^*_{Z(A)}(\Der(A),\Der(A))$$
is a homomorphism for both brackets $[\quad,\quad]^\wedge $ and
$[\quad,\quad]$.
If we denote for the moment also by
$\bar\ze:\Om^*_\Der(A) \to C^*_{Z(A)}(\Der(A),A)$
the embedding of graded differential algebras, then we have
$$
\bar\ze\o i_K = i_{\ze K}\o\bar\ze,\quad \bar\ze\o\L_K = \L_{\ze K}\o\bar\ze:
     \Om^*_\Der(A)\to C^{*}_{Z(A)}(\Der(A),A),
$$
so all the formulas from section \nmb!{5} continue to hold on
$\Der(A,\Om^{*}_\Der(A))$.
\endproclaim

\demo{Proof}
This is obvious from the considerations above.
\qed\enddemo

\subhead\nmb.{8.12} \endsubhead
For $K\in \Der(A,\Om^k(A))$ we consider
$\ze_k\o K\in \Der(A,C^k_{Z(A)}(\Der(A),A))$, and
the corresponding element
$$\ze_k(K)\in C_{Z(A)}(\Der(A),\Der(A))\cong \Der(A,C^k_{Z(A)}(\Der(A),A))$$
then we have:

\proclaim{Lemma}
For $\om\in \Om^q(A)$ and $K\in\Der(A,\Om^{k}(A))$ we have
$$i_{\ze_k(K)}(\ze_q\om) = \ze_{q+k-1}(j_K\om).\tag1$$
\endproclaim

\demo{Proof} Both sides,
$$
i_{\ze_k(K)}\o\ze, \ze\o j_K :\Om^*(A)\to C^{*+k-1}_{Z(A)}(\Der(A),A),
$$
are derivations over $\ze:\Om(A)\to C_{Z(A)}(\Der(A),A)$,
a homomorphism of graded differential algebras,
and both vanish on
$A=\Om^0(A)$, thus it remains to show that they are equal on
$\Om^1(M)$. But for $\om\in\Om^1(A)$ we have
by \nmb!{5.2}.\thetag1
$$\align
(i_{\ze_k(K)}(\ze_1\om))(X_1,\dots,X_k) &=
     (\ze_1\om)(\ze_kK(X_1,\dots,X_k)) \\
&=j_{(\ze_kK(X_1,\dots,X_k))}\om \\
&=j_{(j_{X_k}\dots j_{X_1}K)}\om \\
&=j_{X_k}\dots j_{X_1}j_{K}\om \\
&= \ze_k(j_{K}\om)(X_1,\dots,X_k)\quad \qed
\endalign$$
\enddemo

\proclaim{\nmb.{8.13}. Corollary}
For $K, L\in \Der(A,\Om(A))$ and $\om\in\Om(A)$ we have:
\roster
\item $\ze([K,L]^\De) = [\ze(K),\ze(L)]^\wedge$ for the algebraic
       brackets.
\item $\ze\L_K\om = \L_{\ze(K)}\ze\om$ for the Lie derivations.
\item $\ze([K,L])=[\ze(K),\ze(L)]$ for the Fr\"olicher-Nijenhuis
       brackets.
\endroster
\endproclaim

\head\totoc\nmb0{9}. The differential graded algebra
$\Om_{\text{Out}}(A)$ \endhead

\subhead\nofrills \nmb.{9.1}. The differential graded algebra
$\Om_\Out(A)$ {}\endsubhead
is the subspace of all forms $\om\in \Om_\Der(A)$ which
are basic with respect to to all inner derivations of $A$.
In more detail:

For $a\in A$ let $\ad(a):A\to A$ be given by $\ad(a)b=[a,b]=ab-ba$.
Then the space $\Int(A)$ of all these \idx{\it  inner derivations}
$\ad(a)$ for $a\in A$ is an ideal in $\Der(A)$, and the quotient
$\Out(A)=\Der(A)/\Int(A)$ is called the Lie algebra of outer
derivations.

Then we define $\Om^k_\Out(A)$ to be the set of all
$\om\in \Om^k_\Der(A)$ which satisfy $i_{\ad(a)}\om=0$ and
$\L_{\ad(a)}\om=0$ for all $a\in A$. It is easily seen to be a
differential graded subalgebra.

\proclaim{\nmb.{9.2}. Lemma} The homomorphism
$\ze:\Om_\Der(A)\to C_{Z(A)}(\Der(A),A)$ is part of the following
commutative diagram
$$\CD
\Om_\Der(A) @>\ze>> C_{Z(A)}(\Der(A),A) @>>> C(\Der(A),A)\\
@AAA                @AAA                     @AAA \\
\Om_\Out(A) @>\ze>> C_{Z(A)}(\Out(A),Z(A)) @>>> C(\Out(A),Z(A))
\endCD$$
where in the lower row we have the graded differential subalgebras of
those element which are basic with respect to all inner derivations
of $A$.
\endproclaim

\demo{Proof} Clearly the homomorphisms in the upper row of the
diagram map elements which are basic with respect to all inner
derivations to themselves. It just remains to check that in
$C_{Z(A)}(\Der(A),A)$  and in $C(\Der(A),A)$ these elements form the
sets $C_{Z(A)}(\Out(A),Z(A))$ and $C(\Out(A),Z(A))$, respectively.
Let $\om\in C_{Z(A)}(\Der(A),A)$ or $C(\Der(A),A)$ be basic.
Since $i_{\ad(a)}\om=0$ for all $a\in A$, the skew symmetric
multilinear mapping $\ze(\om):\Der(A)^l\to A$ factors to
$\Out(A)\to A$. By formula \nmb!{5.6} we have
$$\align
0&=(\L_{\ad(a)}\om)(X_1,\dots,X_l) \\
&= \ad(a)(\om(X_1,\dots,X_l) -
\sum\om(X_1,\dots,[\ad(a),X_i],\dots,X_l)\\
&= \ad(a)(\om(X_1,\dots,X_l) - 0
\endalign$$
for all $a\in A$, since $[\ad(a),X_i]\in\Int(A)$. But this means that
$\om$ has values in $Z(A)$.
\qed\enddemo

\proclaim{\nmb.{9.3}. Proposition} The graded differential algebra
$(C_{Z(A)}(\Out(A),Z(A)),d)$ is Morita invariant.
\endproclaim

\demo{Proof}
$Z(A)= HH^1(A;A)$ and $\Out(A)=HH^1(A;A)$, where $HH^*(A;A)$ is the
Hochschild cohomology of $A$ with values in $A$, and also the action
of $\Out(A)$ on $Z(A)$ is described via Hochschild cohomology. Since
Hochschild cohomology is Morita invariant, the result follows.
\qed\enddemo

\enddocument